\documentclass[journal,onecolumn]{IEEEtran}
\usepackage[caption]{subfig}
\usepackage{threeparttable}
\usepackage{graphicx}
\usepackage{amsmath,amssymb,bm,mathrsfs,xfrac,mathtools,mathrsfs,amsthm,amsfonts,bbm}
\usepackage{color, colortbl}
\usepackage[export]{adjustbox}
\usepackage{booktabs} 
\usepackage{tikz} \usetikzlibrary{calc,shapes,arrows,automata,positioning} 
\usepackage{float}
\usepackage{cite}
\usepackage[hyphens]{url}
\usepackage[hidelinks]{hyperref}
\usepackage{soul}
\usepackage{blindtext}
\usepackage{multirow}
\usepackage{pgf}
\usepackage{dashrule}
\usepackage[ruled,vlined]{algorithm2e}
\SetKw{Tri}{$\triangleright$}
\usepackage{enumerate}
\usepackage{tabularx}
\usepackage{array}
\newcolumntype{L}[1]{>{\raggedright\arraybackslash}m{#1}}
\newcolumntype{Y}[1]{>{\raggedright\arraybackslash}m{#1}}
\newcolumntype{R}[1]{>{\vspace{2mm}\raggedright\arraybackslash}m{#1}}
\newtheorem{definition}{Definition}

\newtheorem{proposition}{Proposition}
\newtheorem{lemma}{Lemma}
\newtheorem{example}{Example}

\usepackage{hyperref}
\usepackage{cleveref}
\crefname{lemma}{Lemma}{Lemmas}
\crefname{theorem}{Theorem}{Theorems}
\crefname{algorithm}{Algorithm}{Algorithms}
\crefname{table}{Table}{Tables}
\definecolor{NREL1}{RGB}{000,166,222}

\usepackage{bookmark}

\ifCLASSINFOpdf
\else
\fi

\begin{document}
\bstctlcite{IEEEexample:BSTcontrol}
\title{Unlocking Deep Demand Flexibility\\
via Dynamic Signals}
\author{Xinyang~Zhou, 
        Jing~Shang, 
        Andrey~Bernstein, 
        Stefan~Wager, 
        Moody~Saleh, 
        and~Lara~Pierpoint%
\thanks{Xinyang Zhou is with the Power System Engineering Center, National Laboratory of the Rockies, email: xinyang.zhou@nlr.gov.}%
\thanks{Jing Shang and Stefan Wager are with Operations, Information, and Technology at the Stanford Graduate School of Business, email: \{jshang21, swager\}@stanford.edu.}%
\thanks{Andrey Bernstein is with Eaton Research Labs, email: andreybernstein@eaton.com.}%
\thanks{Moody Saleh and Lara Pierpoint are with Trellis Climate at Prime, email: \{moody.saleh, lara\}@primecoalition.org.}
\thanks{This work was authored by the National Laboratory of the Rockies for the U.S. Department of Energy (DOE), operated under Contract No. DE-AC36-08GO28308. This work was supported by Trellis/Prime Coalition Demand Flexibility project. The views expressed in the article do not necessarily represent the views of the DOE or the U.S. Government. The U.S. Government retains and the publisher, by accepting the article for publication, acknowledges that the U.S. Government retains a nonexclusive, paid-up, irrevocable, worldwide license to publish or reproduce the published form of this work, or allow others to do so, for U.S. Government purposes.}%
}

\maketitle

\begin{abstract}
The rapid proliferation of distributed energy resources (DERs) and the electrification of residential loads offer significant potential for grid flexibility but pose stability challenges under static pricing regimes. Specifically, high levels of automation under static Time-of-Use (TOU) tariffs often induce ``device synchronization,'' where simultaneous responses from home energy management systems (HEMS) create artificial demand peaks that threaten grid stability. This paper proposes a privacy-preserving, one-way dynamic signaling framework to unlock deep demand flexibility from HEMS. We utilize a feedback-based learning algorithm that updates day-ahead price profiles based on aggregate substation demand and environmental contexts, effectively closing the loop between utility objectives and aggregated edge behaviors. The framework is rigorously validated using high-fidelity simulations on an 84-bus distribution network populated with hundreds of HEMS controlling diverse devices, including HVAC, PV, batteries, and flexible loads. Results demonstrate that the proposed mechanism achieves substantial reductions in both peak demand and total load variation. Extensive analyses across diverse climates and scalable deployments confirm the framework's robustness, indicating that dynamic pricing acts as a force multiplier for DERs, with peak shaving potential increasing significantly under high renewable penetration scenarios.
\end{abstract}


\section{Introduction}

The concept of demand response (DR) is undergoing a fundamental paradigm shift. Historically, DR was defined by occasional, discrete events aimed at peak shaving to ensure system reliability. However, the modern grids characterized by the widespread adoption of rooftop photovoltaics (PV), battery energy storage, HVACs, and electric vehicles, requires a broader definition of DR that encompasses continuous, economic, and technical grid management \cite{mathieu2025new}. Recent industry analyses underscore the sheer scale of this latent resource; for example, evaluations of over 25 million residential energy shifts by Renew Home demonstrate that prioritizing subtle, comfort-bounded adjustments can aggregate to over 4 gigawatts of flexible capacity---rivaling the output of major conventional power plants \cite{renewhome2025}. As identified in both academic surveys and industry deployments, there is a critical need for mechanisms that can harness this latent flexibility not just for emergency curtailment, but for efficient daily grid operations.

\subsection{The Challenge of Static Signaling}
A primary barrier to unlocking demand flexibility is the limitation of current signaling mechanisms. Most residential flexibility programs rely on static Time-of-Use (TOU) tariffs \cite{xcel2023residential}. While effective for manual load shifting, TOU structures can become counterproductive in the presence of  automated Home Energy Management Systems (HEMS). When hundreds or thousands of automated devices optimize against the same static price schedule, they often react simultaneously---for instance, pre-cooling homes or charging batteries the moment an off-peak window opens. This phenomenon, known as device synchronization, can create artificial ``timer peaks'' that may exceed the original load, replacing one grid stressor with another.

\subsection{The Communication and Privacy Trade-off}
To mitigate synchronization, theoretical literature often proposes real-time, two-way negotiation markets or direct load control. While optimal in theory, these architectures face significant practical hurdles: they require expensive, low-latency communication infrastructure and raise substantial privacy concerns by exposing detailed household consumption data to utilities. Such complex or intrusive implementation requirements frequently result in low enrollment rates, ultimately rendering the flexibility program less effective as the limited participation fails to reach the critical mass necessary for meaningful grid-level load-shifting impact.

\subsection{Related Work}

The transition from passive distribution networks to active grids populated with DERs has necessitated a fundamental re-evaluation of DR mechanisms. While early DR implementations focused on coarse, event-based peak shedding for reliability, modern frameworks increasingly view flexibility as a continuous grid asset capable of providing ancillary services, congestion management, and renewable integration \cite{mathieu2025new, samad2016automated, vardakas2015survey}.

\subsubsection{The Synchronization Challenge in Pricing Mechanisms}
Static TOU tariffs remain the most prevalent mechanism for residential DR due to their simplicity and regulatory familiarity. However, the rigidity of static pricing creates significant operational risks in high-penetration scenarios where large populations of automated devices optimizing against identical price signals tend to synchronize their consumption, creating ``timer peaks'' that can exceed the original load.

\subsubsection{Centralized vs. Decentralized Architectures}
To address the limitations of static pricing, the literature can primarily be classified into centralized direct control and decentralized market-based coordination \cite{kok2016society}.

\textit{Centralized Control and Aggregation:} Centralized approaches, such as direct load control or model-based aggregation, offer high theoretical performance by treating the grid as a controllable entity. Parvania et al. \cite{parvania2013optimal} demonstrated effective DR aggregation in wholesale markets, while Granitsas et al. \cite{granitsas2025controlling} recently showcased the use of air conditioners for frequency regulation. However, these methods face scalability barriers. They typically require full network observability \cite{nazir2022grid}, reliable low-latency communication infrastructure \cite{zandi2018home}, and sufficient enrollment of customers, all of which act as barriers to entry for distribution utilities.

\textit{Decentralized and Transactive Systems:} Conversely, decentralized approaches such as Peer-to-Peer trading \cite{morstyn2019bilateral, kim2020p2p} and transactive energy systems move decision-making to the edge. While these architectures improve privacy and autonomy, they introduce significant computational overhead and require complex market clearing mechanisms that are difficult to implement within current regulatory frameworks.

\subsubsection{The Algorithmic Gap}
Bridging these architectures requires efficient control algorithms. Extensive surveys on distributed optimization \cite{molzahn2017survey} highlight the trade-off between communication overhead and optimality. Progress has been made in online optimization for networked DERs; for example, \cite{bernstein2019,zhou2019online} developed frameworks that handle time-varying objectives and constraints in real-time. However, while these approaches offer robust theoretical guarantees, there remains a gap for lightweight, privacy-preserving mechanisms that achieve load-shaping benefits without the infrastructure costs of real-time two-way communication. Our work leverages recent advances in feedback-based mechanism design to fill this gap, utilizing one-way dynamic signals to implicitly coordinate massive populations of DERs.

\subsection{Proposed Framework and Contributions}
This paper addresses these challenges by developing and validating a privacy-preserving, one-way dynamic signaling framework. Unlike static tariffs, our approach utilizes a feedback-based dynamic algorithm that updates 24-hour ahead price profiles based on aggregate substation demand and weather contexts. This allows the utility to ``shape'' demand dynamically, smoothing out volatility and preventing synchronization without two-way communication or requiring direct control of customer devices.

This work demonstrates a rigorous, data-driven validation of this framework within a high-fidelity simulation environment. Moving beyond stylized models, we employ a realistic distribution system (SMART-DS) \cite{nrel2023smartds} and historical weather data \cite{zippenfenig2023open} to quantify performance. Key contributions include:

\begin{itemize}
    \item \textbf{Empirical Validation of Efficacy:} We demonstrate that the proposed one-way signaling mechanism achieves performance metrics remarkably close to ideal two-way control scenarios, delivering a 16.01\% Peak Demand Shaving (PDS) and 19.33\% variation reduction in summer months.
    \item \textbf{Robustness Across Climates:} Sensitivity analysis across three distinct climatic zones (Denver, Los Angeles, and Phoenix) confirms the framework's effectiveness under varying cooling loads, with variation reductions reaching up to 40\% in extreme climates like Phoenix.
    \item \textbf{Scalability Analysis:} We verify the framework's applicability to large-scale systems, showing similar performance when scaling from hundreds to over 14,000 participating households.
    \item \textbf{Synergy with Renewables:} We show that dynamic pricing acts as a force multiplier for DERs; increasing PV and battery penetration to 60\% nearly doubles the peak shaving capability to 32.79\%.
\end{itemize}

The remainder of this paper is organized as follows: Section II formulates the system model and feedback framework. Section III briefly outlines the signal computation algorithm. Section IV presents extensive numerical results and sensitivity analyses. Section V concludes this work.

\section{System Modeling And Problem Formulation}
This section establishes the mathematical modeling for the proposed demand flexibility framework. We define a comprehensive system that models the interaction between a utility and a population of residential consumers, framed as a Stackelberg game bi-level optimization problem. The objective of the framework is to optimize overall grid performance—--specifically by smoothing the demand profile and reducing peak loads--—while respecting individual consumer preferences and the physical constraints of behind-the-meter devices. To this end, we provide detailed physical models for various smart devices, including HVAC systems, flexible loads, battery storage, and rooftop PV, and characterize the collective properties of the household population.

\subsection{System Overview} \label{sec:schemeOverview}
The architecture of our framework is centered on a closed-loop, one-way signaling mechanism, as illustrated in Fig. 1. The operational flow proceeds through the following daily cycle:
\subsubsection{Signal Generation}  At the beginning of each day, the utility's price signal generator utilizes exogenous inputs—such as weather and solar irradiance forecasts—alongside the demand profile from the previous 24 hours to generate an optimized 24-hour price signal vector.
\subsubsection{One-Way Communication} These price signals are broadcast through a one-way communication channel to all households within the distribution network. This approach is chosen for its practical feasibility and its ability to preserve consumer privacy by avoiding the need for complex, iterative negotiations.
\subsubsection{HEMS Optimization} We assume that each modeled household is equipped with a HEMS. Participating households utilize the received signals to solve a local optimization problem, balancing electricity costs against user comfort. Non-participating households  continue to operate solely based on their internal preferences, regardless of the utility's price signal.
\subsubsection{Physical System Feedback} The individual power consumption profiles are aggregated by the physical power system, resulting in a demand profile observed at the substation. This realized demand then serves as feedback to the utility, informing the signal generation process for the subsequent 24-hour horizon and effectively closing the control loop.

The subsequent subsections provide a formal mathematical characterization of each system component.

\begin{figure}[!h]
\centerline{\includegraphics[width=.6\textwidth]{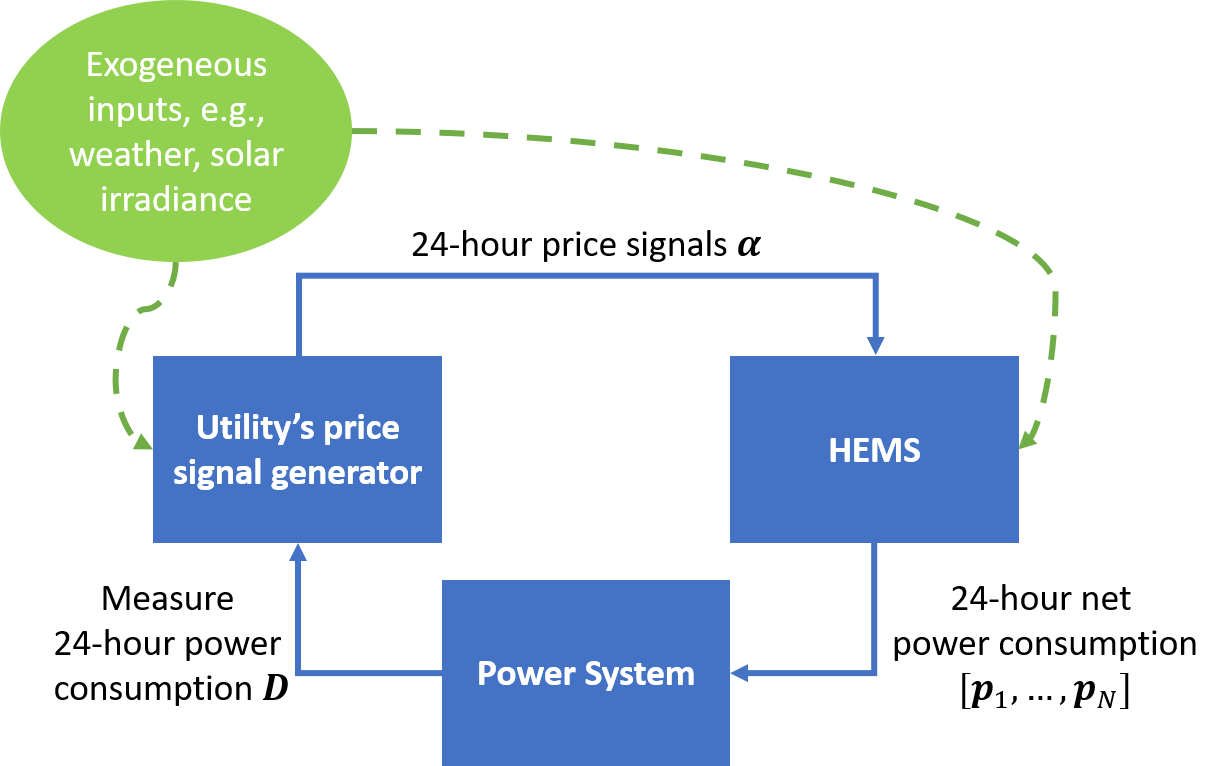}}
\caption{The system overview: Utility generates day-ahead 24-hour price signals $\bm{\alpha}$ based on previous day's 24-hour demand profile, as well as exogenous inputs including weather forecast, and sends them to HEMS. HEMS calculates optimized power consumption given $\bm{\alpha}$, its own preferences, and exogenous inputs. The physical power system itself helps close the loop by providing feedback to utility for the next day's price decision.}
\label{fig:flowchart}
\end{figure}

\subsection{Power System Modeling}
In this work, we consider a distribution network populated with $N$ households. We collect and index all $N$ households in the system by a set $\mathcal{N}:= \{1, \ldots, N\}$. We assume that households have the freedom to opt in and out of the demand response program by choosing to integrate price signals into their control decisions or not. We denote by $\mathcal{N}_{in}\subseteq\mathcal{N}$ and $\mathcal{N}_{out}\subseteq\mathcal{N}$ the subsets of households that opt in and out of the program, respectively. We also assume $\mathcal{N}_{in}$ and $\mathcal{N}_{out}$ are non-overlapping, and that $\mathcal{N}_{in}\bigcup\mathcal{N}_{out}=\mathcal{N}$. Since we assume one household is equipped with one HEMS, we will use the terms household and HEMS interchangeably for convenience.

Denote by $D\in \mathbb{R}$ the total demand at the substation. We divide a day equally into $T$ slots, with each slot lasting $\Delta$, e.g., 1 hour. We let $\bm{D} = [D^1, \ldots, D^t, \ldots, D^T]^\top$ denote a profile collecting total demand for all time slots of the day.

We next formulate the general demand flexibility model with two  players: the utility and households. The problem is formulated as a Stackelberg game bi-level optimization of the form:
\begin{align}
    &\min_{\bm{\alpha}} \rho\left(\bm{D}(\bm{\alpha})\right). \label{eq:sgame}
\end{align} 
Here $\bm{\alpha}=[\alpha^1, \ldots, \alpha^t, \ldots, \alpha^T]^{\top}\in\mathbb{R}^{T}$ is the vector collecting price signals for all the time slots, with each $\alpha^t$ being either actual or virtual price for time slot $t$ as explained below; the objecitve function $\rho (\cdot)$ measures the quality of the overall demand profile based on utility needs; and $\bm{D}(\bm{\alpha})$ is the demand profile realized by the households for a given signal $\bm{\alpha}$. As explained below, $\bm{D}(\bm{\alpha})$ is itself an outcome of a HEMS optimization problem, resulting in a bi-level optimization formulation. See Fig.~\ref{fig:flowchart} for an illustration of the system overview.

\subsection{Utility's Objective}\label{utility_objective}
When employing a demand flexibility program, a utility may have multiple objectives. In our framework, we express these objectives via a single objective function term $\rho(\bm{D})$, which can be configured based on specific utility needs. In our experiments, we used the following objective function that reflects the goal of 
\emph{lower peak hour demand while smoothing out the entire demand profile}:
\begin{eqnarray}\label{eq:costfunction}
    &\rho(\bm{D}) = \epsilon \cdot QV(\bm{D}) + (1-\epsilon)\cdot\|\bm{D}\|^2_2,
\end{eqnarray}
where the quadratic variation $QV(\bm{D})$ is defined as the squared temporal demand difference from one hour to the next:
\begin{eqnarray}
    QV(\bm{D}):=\sum_{t=1}^{T-1}(D^{t+1}-D^{t})^2.
\end{eqnarray}
Here, $\epsilon\geq 0$ is a chosen parameter to balance the two objectives. This specific formulation is selected for its smoothness and strong convexity, which enhances numerical performance during optimization, though the framework allows for other similar objective functions. To practically realize these objectives within the distribution network, we will introduce an algorithm in Section~\ref{sec:alg} that delivers these grid-level goals through the optimized design of dynamic pricing signals.

\begin{figure}[!h]
\centerline{\includegraphics[width=.9\textwidth]{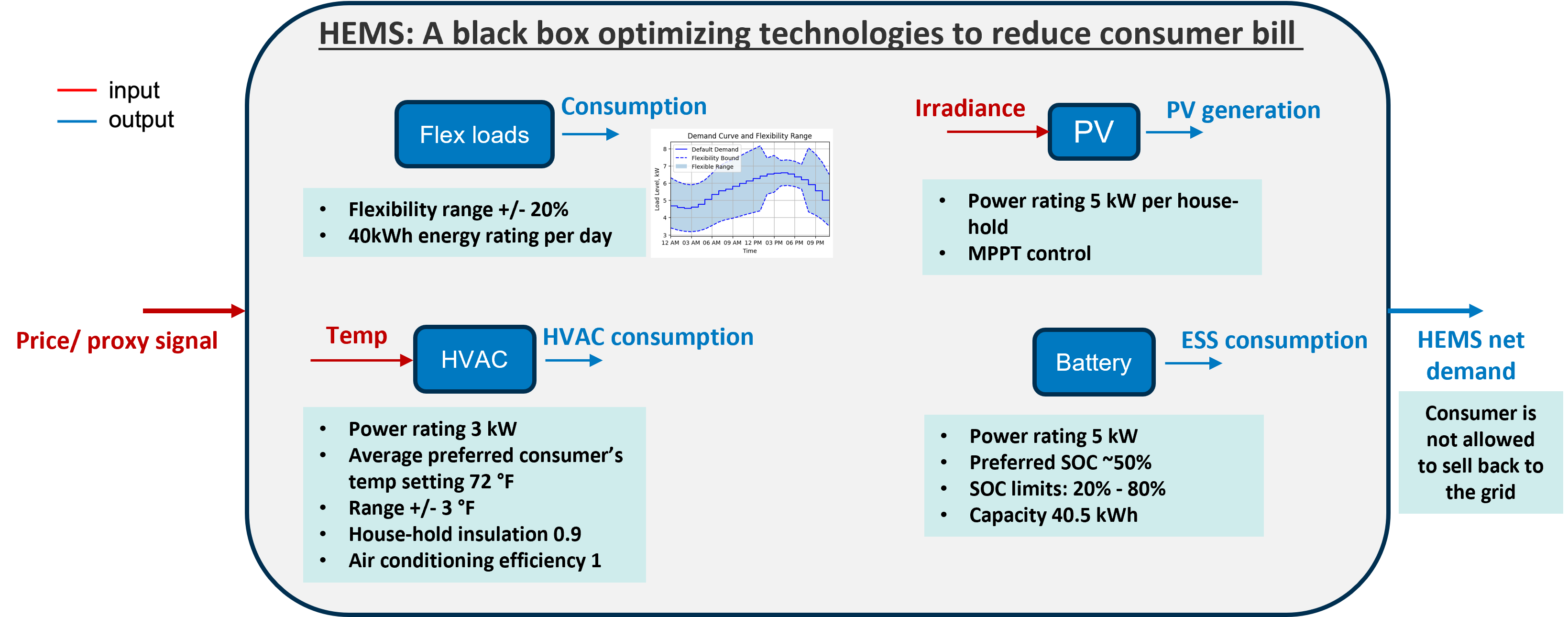}}
\caption{HEMS optimizes power consumption profiles of its devices given price signals.}
\label{fig:HEMS}
\end{figure}

\subsection{Household Modeling}

We assume that a household $i\in\mathcal{N}$ has a few controllable devices indexed by $d$, collected within a device set $\mathcal{D}_i$. Device $d\in\mathcal{D}_i$ consumes power $\bm{p}_{i,d}=[p^1_{i,d}, \ldots, p^t_{i,d}, \ldots, p^T_{i,d}]^\top\in\mathbb{R}^{T}$ over a day, adding up to a household consumption profile of
\begin{eqnarray}
\bm{p}_i=\sum_{d\in\mathcal{D}_i}\bm{p}_{i,d}.
\end{eqnarray}

Each device's consumption profile $\bm{p}_{i,d}$ should fall within the ranges of physical limits and user preference characterized by a convex feasible set $\Omega_{i,d}$.
We define a cost function $C_{i,d}(\bm{p}_{i,d})$ for each device and an aggregated cost function for the entire household as
\begin{eqnarray}
\sum_{d\in\mathcal{D}_i} C_{i,d}(\bm{p}_{i,d}).
\end{eqnarray}
The specific forms of $\Omega_{i,d}$ and $C_{i,d}$ will be detailed when we model individual devices in the next subsection.

\subsubsection{Participating Households}
Given a signal vector $\bm{\alpha}$, a participating household $i\in\mathcal{N}_{in}$ optimizes its power consumption plan $\bm{p}_i$ by solving the following optimization problem to balance between its own power usage preferences $\sum_{d\in\mathcal{D}_i} C_{i,d}(\bm{p}_{i,d})$ and signal responsiveness (e.g., electricity bill) $\bm{\alpha}^{\top} \bm{p}_i$:
\begin{subequations}\label{eq:opt1}
\begin{eqnarray}
        \min_{\bm{p}_{i,d}, d\in\mathcal{D}_i} && \sum_{d\in\mathcal{D}_i} C_{i,d}(\bm{p}_{i,d}) + \bm{\alpha}^{\top} \bm{p}_i \label{eq:hemscost1}\\
        \text{s.t.}&& \bm{p}_{i,d}\in\Omega_{i,d}, d\in\mathcal{D}_i,\\
        && \bm{p}_i=\sum_{d\in\mathcal{D}_i}\bm{p}_{i,d},\\
        && g(\bm{p}_i)\leq 0,
\end{eqnarray}
\end{subequations}
where $g(\bm{p}_i)$ is a convex function representing some optional household-level constraints, e.g., total power rating limits.

\subsubsection{Non-Participating Households}
For a non-participating household $i\in\mathcal{N}_{in}$, it instead solves the following optimization problem regardless of the price signals:
\begin{subequations}\label{eq:opt2}
\begin{eqnarray}
        \min_{\bm{p}_{i,d}, d\in\mathcal{D}_i} && \sum_{d\in\mathcal{D}_i} C_{i,d}(\bm{p}_{i,d}) \label{eq:hemscost1}\\
        \text{s.t.}&& \bm{p}_{i,d}\in\Omega_{i,d}, d\in\mathcal{D}_i,\\
        && \bm{p}_i=\sum_{d\in\mathcal{D}_i}\bm{p}_{i,d},\\
        && g(\bm{p}_i)\leq 0,
\end{eqnarray}
\end{subequations}
where the difference is the lack of response towards signal in the cost function.

We show in Fig.~\ref{fig:HEMS} how a HEMS coordinates multiple devices given control signals and user preferences. We next define the mathematical model of individual devices for HEMS to optimize and control.

\subsection{Smart Devices Modeling}\label{sec:devicemodel}

We now define the models of smart devices under control, i.e., the specific mathematical forms of $C_{i,d}$ and $\Omega_{i,d}$.
\subsubsection{HVAC}
For an HVAC $d\in\mathcal{D}_i$, let $\bm{T}^{out}\in\mathbb{R}^{T}$ denote the outside temperature forecast for the next $T$ time slots, and $\bm{T}_{i,d}^{in}\in\mathbb{R}^{T}$ the room temperature to be monitored and controlled by HVAC power consumption $\bm{p}_{i,d}$. The room temperature dynamics from time $t$ to $t+m$ can be modeled in the following form \cite{li2011optimal}:
\begin{eqnarray}
T_{i,d}^{in, t+m}&\hspace{-7pt}=&\hspace{-7pt}T_{0,i,d}^{t+m}+\sum_{\tau=0}^{m-1}(1-\zeta_1)^{m-1-\tau}
\zeta_2 p_{i,d}^{t+\tau}, \label{eq:tp1}
\end{eqnarray}
with $\zeta_1$ and $\zeta_2$ are parameters specifying the HVAC cooling effectiveness and the thermal characteristics of the house insulation, respectively. 
Here, the parameter $T_{0,i,d}^{t+m}$ characterizes the natural indoor temperature if HVAC is not on:
\begin{eqnarray}
\hspace{-3mm}T_{0,i,d}^{t+m}&\hspace{-7pt}=&\hspace{-7pt}(1-\zeta_1)^m T_{i,d}^{in,t}+\sum_{\tau=0}^{m-1}(1-\zeta_1)^{m-1-\tau} \zeta_2 T^{t+\tau}_{out,i,d}. \label{eq:tp2}
\end{eqnarray}
 We assume a preferred temperature vector $\bm{T}^{prefer}_{i,d}\in\mathbb{R}^T$ and a cost function penalizing indoor temperature deviating from it with a customized elasticity parameter $\gamma_{i,d}$ as:
 \begin{eqnarray}
 C_{i,d}(\bm{p}_{i,d}) = \gamma_{i,d}\|\bm{T}^{in}_{i,d}(\bm{p}_{i,d}) - \bm{T}^{prefer}_{i,d} \|_2^2, \label{eq:cost_hvac}
\end{eqnarray}
We also set hard temperature constraints with customized upper and lower temperature bounds $\overline{\bm{T}}_{i,d}$ and $\underline{\bm{T}}_{i,d}$:
\begin{eqnarray}
\underline{\bm{T}}_{i,d}\leq \bm{T}_{i,d}^{in}\leq \overline{\bm{T}}_{i,d}.\label{eq:tp3}
\end{eqnarray}
Here, $\gamma_{i,d}$ characterizes how important household $i$ views the temperature preferences, in comparison to signal responsiveness. 

We  assume that the HVAC's power consumption is bounded by the maximum power rating $\overline{\bm{p}}_{i,d}$ as:
\begin{eqnarray}
\bm{0}\leq \bm{p}_{i,d}\leq \overline{\bm{p}}_{i,d}.\label{eq:hvac_power}
\end{eqnarray}

The feasible set of an HVAC $i, d$ can therefore be defined as:
\begin{align} 
\Omega_{i,d} = \left\{ \bm{p}_{i,d} \big|  \eqref{eq:tp1}-\eqref{eq:tp2}, \eqref{eq:tp3}-\eqref{eq:hvac_power} \right\}.\label{eq:const_hvac}
\end{align}

\subsubsection{Flexible Load}
For a flexible load $d\in\mathcal{D}_i$, we assume there to be a preferred power consumption schedule $\bm{p}_{i,d}^{prefer}$ with a fixed total energy consumption $E_{i,d}=\sum_{t=1}^T p_{i,d}^{prefer,t}$ within the decision horizon $T$. We assume the flexible load can be increased and decreased within certain bounds for each hour as\footnote{Note that the household's non-controllable power consumption can all be aggregated into the flexible load model by adding to the lower bound $\underline{\bm{p}}_{i,d}$.}:
\begin{eqnarray}
\underline{\bm{p}}_{i,d} \leq \bm{p}_{i,d}\leq \overline{\bm{p}}_{i,d},\label{eq:load_power}
\end{eqnarray}
with a cost function penalizing power deviation from the preferred schedule:
 \begin{eqnarray}
 C_{i,d}(\bm{p}_{i,d}) = \gamma_{i,d}\|\bm{p}_{i,d} - \bm{p}^{prefer}_{i,d} \|_2^2, \label{eq:cost_load}
\end{eqnarray}
while maintaining the total energy consumption throughout the scheduled time horizon:
\begin{eqnarray}
\sum_{t=1}^T p_{i,d}^{t} = E_{i,d}.\label{eq:totalload}
\end{eqnarray}
As a result, we do not reduce daily energy consumption from flexible load, but energy usage shift. Similarly, $\gamma_{i,d}$ is an elasticity parameter indicating household's trade-off preference between its preferred schedule consumption plan and utility signals responsiveness.
We can summarize the feasible set of the flexible load $i, d$ as 
\begin{align} 
\Omega_{i,d} = \left\{ \bm{p}_{i,d} \big|  \eqref{eq:load_power}, \eqref{eq:totalload}\right\}.\label{eq:const_load}
\end{align}

\subsubsection{Battery}
For a {battery} 
$d\in\mathcal{D}_{i}$, let $\bm{SOC}_{i,d}\in\mathbb{R}_+^T$ represent its non-negative state of charge (SOC) profile over the day. The dynamics for SOC evolution from time $t$ to $t+m$ take the form\footnote{Here we have assumed that the charging and discharging efficiencies are 1. 
We can deal with arbitrary battery efficiencies in the same way as in \cite{stai2017dispatching} that approximates the losses using an extended model; however, for simplicity of exposition, we consider an efficiency of 1.} of:
\begin{eqnarray}
SOC_{i,d}^{t+m}=SOC_{i,d}^t+\sum_{\tau=0}^{m-1} p_{i,d}^{t+\tau}, \label{eq:soc1}
\end{eqnarray}
with power limits within maximum charging rating $\overline{\bm{p}}_{i,d}\geq \bm{0}$ and maximum discharging rating $\underline{\bm{p}}_{i,d}\leq \bm{0}$:
\begin{eqnarray}
\underline{\bm{p}}_{i,d}\leq \bm{p}_{i,d}\leq \overline{\bm{p}}_{i,d}.\label{eq:soc2}
\end{eqnarray}
Moreover, $\bm{SOC}_{i,d}$ should fall within an acceptable SOC bounds of $\underline{\bm{SOC}}_{i,d}$ and $\overline{\bm{SOC}}_{i,d}$ (e.g., 20\% and 80\%) as:
\begin{eqnarray}
\underline{\bm{SOC}}_{i,d}\leq \bm{SOC}_{i,d}\leq \overline{\bm{SOC}}_{i,d}.\label{eq:soc3}
\end{eqnarray}

Here, the constraints can involve charging deadline, e.g.,
charge the EV battery with at least 90\% SOC by 8 a.m. We assign a cost function to penalize deviations of SOC levels from preferred levels $\bm{SOC}_{i,d}^{prefer}$ with a chosen elasticity parameter $\gamma_{i,d}$ for this battery as:
\begin{eqnarray}
 C_{i,d}(\bm{p}_{i,d}) = \gamma_{i,d}\|\bm{SOC}_{i,d}(\bm{p}_{i,d}) - \bm{SOC}^{prefer}_{i,d} \|_2^2. \label{eq:cost_bess}
\end{eqnarray}
We denote feasible set of battery $i, d$  as
\begin{align} 
\Omega_{i,d} = \left\{ \bm{p}_{i,d} \big|  \eqref{eq:soc1}- \eqref{eq:soc3}\right\}.\label{eq:const_bess}
\end{align}

\subsubsection{Roof-Top PVs}
For a PV $d\in\mathcal{D}_i$, let $\bm{I}\in\mathbb{R}^{T}$ be the solar irradiance forecast profile, and the corresponding maximum PV panel generation $\overline{\bm{p}}_{i,d}\in\mathbb{R}^{T}$ for the next day. Note that we assign negative sign to power generation, in contrast to positive sign to power consumption. The feasible region of a PV inverter thus has the following form:
\begin{align}
\Omega_{i,d} =  \left\{ \bm{p}_{i,d} \big| \overline{\bm{p}}_{i,d} \leq \bm{p}_{i,d}  \leq  \bm{0}  \right\}.\label{eq:const_pv}
\end{align}

We assume a convex cost function for PV to penalize power curtailment with a chosen elasticity parameter $\gamma_{i,d}$ for this PV:
\begin{eqnarray}
 C_{i,d}(\bm{p}_{i,d}) = \gamma_{i,d}\|\bm{p}_{i,d} - \overline{\bm{p}}_{i,d} \|_2^2. \label{eq:cost_pv}
\end{eqnarray}




\subsection{One-Way Communication and Control Signals}

Theoretically, in order to find the optimal price via \eqref{eq:sgame}, the utility needs to know the HEMS demand $\bm{D}(\bm{\alpha})$ in response to the signal $\bm{\alpha}$. In the most naive implementation, this requires a two-way communication system between the utility and households. For example, the utility can find an optimal $\bm{\alpha}$ by ``negotiation'' communication rounds, where in the beginning of the day, it iterates with the HEMS and updates $\bm{\alpha}$ opposed to the gradient of the objective function until convergence. However, this approach would be infeasible in practice as it relies on a complex communication infrastructure and requires reliable responses from HEMS in order to converge to an optimal signal. 

Instead, in this work, we assume a one-way communication channel from the utility to households for practicality and preservation of household privacy. 
Particularly, utility utilizes information it has gathered, e.g., past demand records, weather forecast, time of the year, to generate control signals $\bm{\alpha}=[\alpha^1, \ldots, \alpha^T]^{\top}\in\mathbb{R}^{T}$ for the next $T$ time slots, and broadcast them to all households. Participating households in the program will take these signals into consideration, together with  their own preferences and constraints when planning power consumption; non-participating households make decisions only based on their own preferences and constraints regardless of price signals. At the end of the day, the utility observes the resulting demand vector $\bm{D}(\bm{\alpha})$ and updates the signals for the next day. The details of how the signals are computed appear in Section~\ref{sec:alg}.

\subsection{Household Population Properties}
\subsubsection{Participation Rate} 
At the beginning of this section, we have defined $\mathcal{N}_{in} \subseteq \mathcal{N}$ and $\mathcal{N}_{out} \subseteq \mathcal{N}$ as the subsets of households that respectively opt in and out of the demand response program. To quantify the scale of engagement, we define the participation rate, $\mathcal{P}$, as the ratio of participating households to the total population:
\begin{eqnarray}
\mathcal{P} = \frac{|\mathcal{N}_{in}|}{|\mathcal{N}|}.
\end{eqnarray}
Notably, only households in $\mathcal{N}_{in}$ are sensitive to utility signals. Consequently, for a fixed population size $|\mathcal{N}|$, a higher participation rate correlates with greater aggregate responsiveness to signal fluctuations. The impact of varying participation rates on overall demand response performance is further illustrated in Section~\ref{sec:numerical}.

\subsubsection{Elasticity} 
As defined in the cost function \eqref{eq:hemscost1}, each participating household $i \in \mathcal{N}_{in}$ optimizes a trade-off between device-level discomfort/disutility $C_i(\bm{p}_i)$ and total energy expenditures $\bm{\alpha}^{\top} \bm{p}_i$. The sensitivity of this trade-off is governed by the parameters $\gamma_{i,d}$ assigned to individual devices $d$ within their respective cost functions \eqref{eq:cost_hvac}, \eqref{eq:cost_load}, \eqref{eq:cost_bess}, and \eqref{eq:cost_pv}. By aggregating these individual device cost functions, household $i$ exhibits a composite elasticity profile. Furthermore, when these behaviors are summed across the entire participating population, the aggregate group demonstrates collective elasticity characteristics. These emergent population-level properties and their numerical demonstrations are detailed in Section~\ref{sec:numerical}.

\subsubsection{Flexibility} We define flexibility as the physical capacity to reduce or shift energy consumption across different time intervals. Within our modeling framework, this property is fundamentally constrained by the device-level operational requirements—specifically the hard constraints defined in equations \eqref{eq:const_hvac}, \eqref{eq:const_load}, \eqref{eq:const_bess}, and \eqref{eq:const_pv}—which dictate the physical and user preference boundaries of energy usage for each individual household. By aggregating these constraints across all participating households, we characterize the collective flexibility of the entire population. This population-level flexibility is further examined through numerical illustrations in Section~\ref{sec:numerical}.

Note that although households in both $\mathcal{N}_{in}$ and $\mathcal{N}_{out}$ possess innate physical flexibility according to our modeling framework, non-participating households remain at their baseline preferred states by design. Specifically, they maintain $\bm{T}^{prefer}_{i,d}$ for HVAC, $\bm{p}^{prefer}_{i,d}$ for flexible loads, $\bm{SOC}^{prefer}_{i,d}$ for battery storage, and $\overline{\bm{p}}_{i,d}$ for PV generation. Consequently, these households do not utilize their flexibility ranges. In contrast, only participating households adjust their consumption in response to utility signals, thereby actively demonstrating and utilizing their flexibility potential.

\begin{figure}[!h]
\centerline{\includegraphics[width=.4\textwidth]{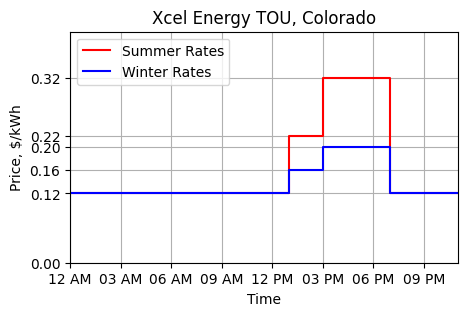}
\includegraphics[width=.4\textwidth]{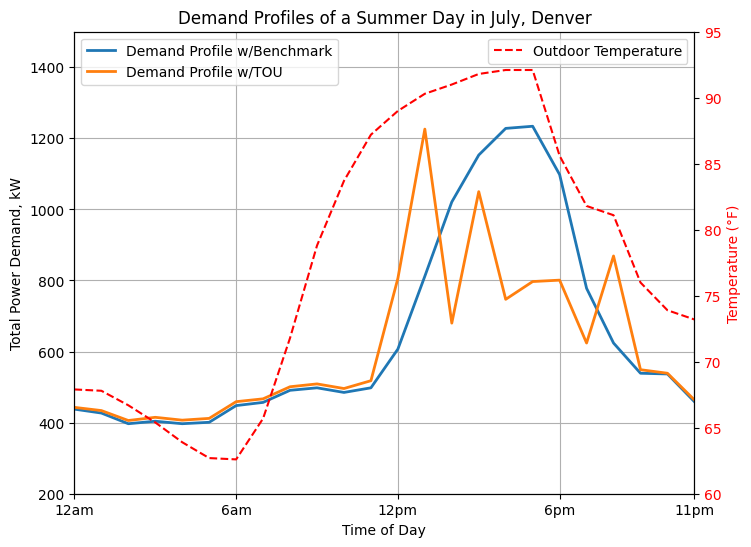}
}
\caption{(Left) TOU used by Xcel Energy in Colorado in 2023; (right) simulated demand profile under TOU in a summer day.}
\label{fig:xceltou}
\end{figure}

\subsubsection{Example: Over Reaction from Mass Population under TOU}
A popular choice for utility pricing is TOU tariff, such as the structure implemented by Xcel Energy in Colorado. With this strategy, electricity rates are pre-set for specific blocks of time—typically categorized as on-peak, mid-peak, and off-peak (see, e.g., Fig.~\ref{fig:xceltou}~(left)). While these static structures are designed to encourage shifting demand away from peak periods, they can become counterproductive in the presence of high-level behind-the-meter automation.

A significant risk is device synchronization: when a large population of households optimize consumption against the same static price signal, the majority of consumers may react simultaneously. For instance, at the precise moment a lower-priced off-peak period begins, HEMS across the network may concurrently trigger high-load devices, such as HVAC systems and battery storage systems. This collective, synchronized reaction can inadvertently create new, secondary demand peaks---sometimes exceeding the original ones---thereby compromising grid stability. Fig.~\ref{fig:xceltou} (right) illustrates such a simulated scenario, where TOU-driven synchronization shifts the load in a poorly contained manner, replacing the original peak with a new, artificial one.

This phenomenon underscores the limitations of static tariffs and motivates the development of the more carefully constructed, dynamic pricing strategy presented in the following section.

\section{Signal Computation Algorithm}\label{sec:alg}

In this section, we present the algorithmic framework for  signal generation and describe how system-level feedback is incorporated to update the pricing policy. This corresponds to Steps 1) and 4) in the system overview outlined in Section \ref{sec:schemeOverview}. We demonstrate two learning algorithms for the optimal signal computation. We first consider a baseline feedback-based learning algorithm introduced in \cite{mehrabi2024optimalmechanismsdemandresponse}. This approach relies solely on aggregated demand observations and is provably mean-field optimal in the special case where the household cost function $C_i(\bm{p}_i)$ is linear in the price vector $\bm{p}_i$. We then extend this framework to a cluster-feedback-based learning algorithm that incorporates environmental and contextual factors into the pricing rule. By leveraging exogenous signals, this approach learns a context-dependent price mapping that adapts to heterogeneous operating conditions.

For both algorithms, we outline the practical implementation steps and discuss their expected performance characteristics as well as potential limitations. Detailed technical specifications—including objective formulations, constraints, training procedures, and hyperparameter configurations—are deferred to the appendix.

\subsection{Feedback-Based Context-Agnostic Algorithm}

In this section, we describe the baseline feedback-based algorithm proposed in \cite{mehrabi2024optimalmechanismsdemandresponse}, which constructs price signals without conditioning on day-specific contextual information. The method operates in an online learning framework, updating the price vector on a daily basis in response to the realized aggregate demand profile. At a high level, the update rule adjusts relative prices across hours to counteract observed demand imbalances: for example, if a demand peak occurs at hour $j$ on a given day, the algorithm increases the relative price at hour $j$ in subsequent iterations to incentivize load shifting and smooth the demand trajectory, thereby mitigating grid risk.
\cref{alg1} presents the main procedure of this feedback-based approach, which relies exclusively on aggregate demand feedback. \cite{mehrabi2024optimalmechanismsdemandresponse} showed that this algorithm recovers mean-field optimal prices wherenever the household cost function $C_i(\bm{p}_i)$ is linear in power $\bm{p}_i$.

\begin{center}
\begin{minipage}{0.75\textwidth}
\begin{algorithm}[H]
\DontPrintSemicolon
\setlength{\algowidth}{0.75\textwidth} 
\newcommand{\rcomment}[1]{\hfill\parbox[t]{3cm}{\raggedleft\footnotesize\textit{#1}}}

\caption{Feedback-Based Context-Agnostic Algorithm} \label{alg1}

\vspace{3mm}
\Tri Initialize price signal $\bm{\alpha}_0$\;
\vspace{3mm}

\For{$t=1,2,\cdots T$}{
\vspace{3mm}
\Tri Send price signal $\bm{\alpha}_t$ to a set of HEMS\;

\Tri Observe the aggregate consumption $\bar{\bm{p}}_t = \frac{1}{n} \sum_{i=1}^n
\bm{p}_i(\bm{\alpha}_t)$\;

\Tri Form a gradient $\bm{g}_t = \bar{\bm{p}}_t $, and use it to update $\bm{\alpha}_{t+1}$\;
}

\vspace{3mm}
\Return Converged optimal signal $\bm{\alpha}_T$

\vspace{3mm}

\end{algorithm}
\end{minipage}
\end{center}

\subsection{Cluster-Feedback-Based Context-Enriched Algorithm}

In this section, we introduce a contextual extension of \cref{alg1} that incorporates exogenous information into the construction of the price signal. While the baseline feedback-based algorithm leverages time-varying electricity prices to mitigate grid risk, it does not explicitly account for contextual factors that may influence the structure of the optimal pricing policy. Consequently, it is limited in its ability to differentiate between days with distinct underlying system conditions.

In practice, optimal prices evolve in response to a variety of complex factors, including weather forecasts, sunrise and sunset times, seasonal effects, and day-of-week patterns. To address this limitation, we propose a cluster-feedback-based algorithm for contextual energy pricing. By conditioning price updates on structured contextual information, the proposed approach learns differentiated pricing rules that adapt to heterogeneous operating environments.

The first step of the proposed framework is to train a classifier using historical weather data. The classifier jointly learns $K$ cluster centers and a probabilistic mapping from weather features to a distribution over clusters. In particular, for each predicted set of contextual features, the model outputs a probability vector on the $K$-dimensional simplex, representing the likelihood of membership in each cluster. For instance, when $K=2$, the learned clusters may correspond to latent regimes such as ``hot sunny days'' and ``cold cloudy days''. A set of forecast characterized by high temperature and strong solar radiation would then be mapped to a distribution such as $(0.9,0.1)$, indicating a high likelihood of belonging to the first regime.

Conditional on the learned cluster structure, we implement the feedback-based learning algorithm within each cluster. The resulting price signal is formed as a weighted average across clusters, where the weights correspond to the probabilities assigned by the classifier to each cluster based on the forecasted contextual variables. Similarly, gradient updates are scaled by the corresponding cluster weights. The detailed procedure of this modified algorithm is presented in \cref{alg2}. When the number of clusters satisfies $K=1$, \cref{alg2} reduces to the baseline feedback-based algorithm described in \cref{alg1}. 

\begin{figure}[ht]
    \centering
    \includegraphics[width=1\linewidth]{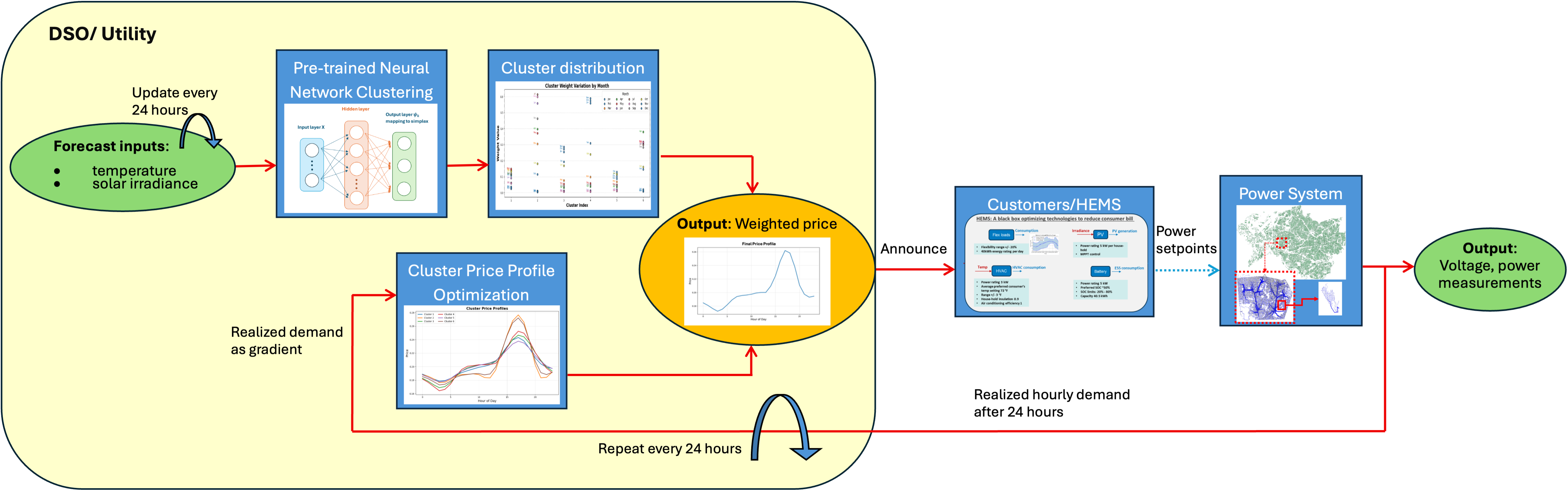}
    \caption{Algorithm Visualization}
    \label{fig:alg_vis}
\end{figure}

\Cref{fig:alg_vis} illustrates the complete pipeline integrating contextual price signal learning with system-level feedback in the power system. A comprehensive description of algorithmic specifications, including parameter configurations and implementation details, is provided in the appendix.

\begin{center}
\begin{minipage}{0.75\textwidth}
\begin{algorithm}[H]
\DontPrintSemicolon
\setlength{\algowidth}{0.75\textwidth} 
\newcommand{\rcomment}[1]{\hfill\parbox[t]{3cm}{\raggedleft\footnotesize\textit{#1}}}

\caption{Cluster-Feedback-Based Context-Enriched Algorithm} \label{alg2}

\vspace{3mm}
\Tri Train a classifier $\psi(\cdot)$ that forms $K$ clusters on a set of historical weather data $X_t$

\vspace{3mm}
\Tri Map contextual information $X_t$ to the distribution on clusters $(\psi(X_t)_1, \cdots \psi(X_t)_K)$, where $\sum_{k=1}^K \psi(X_t)_k = 1$\;
\vspace{3mm}

\Tri Initialize price signal $\bm{\alpha}_{0,1},\cdots,\bm{\alpha}_{0,K}$ for each cluster\;
\vspace{3mm}

\For{$t=1,2,\cdots T$}{
\vspace{3mm}
\Tri Send cluster-weighted price signal $\bm{\alpha}_t:=\sum_{k=1}^K\psi(X_t)_k\bm{\alpha}_{t,k}$ to a set of HEMS\;

\Tri Observe the aggregate consumption $\bar{\bm{p}}_t = \frac{1}{n} \sum_{i=1}^n
\bm{p}_i(\bm{\alpha}_t)$\;

\Tri Form a gradient $\bm{g}_t = \bar{\bm{p}}_t $\;
\For{$k=1,\cdots K$}{
\vspace{1mm}
\Tri Update $\bm{\alpha}_{t,k}$ with cluster-weighted gradient $\psi(X_t)_k\cdot \bm{g}_t$
}
}

\vspace{3mm}
\Return Converged optimal signal for each cluster $\bm{\alpha}_{T,1},\cdots,\bm{\alpha}_{T,K}$

\vspace{3mm}

\end{algorithm}
\end{minipage}
\end{center}

\subsection{Expected Outcome and Challenges}
We conclude this section with a brief discussion of the anticipated performance characteristics of the proposed algorithm, highlighting both its advantages relative to homogeneous pricing frameworks and the conditions under which its benefits may be lessened.

As discussed earlier, the primary motivation of the cluster-feedback-based approach is to account for day-to-day variation in the optimal price signal driven by exogenous factors. Accordingly, we expect the proposed method to outperform algorithms that assume homogeneous day types in environments where contextual heterogeneity substantially affects demand patterns and system conditions. Such scenarios include, but are not limited to: 
\begin{enumerate}[(i)]
    \item Regions with pronounced day-to-day variability in weather conditions, such as substantial fluctuations in temperature or solar irradiance, which induce significant variation in load and renewable generation profiles.
    \item Systems characterized by high penetration of HVAC loads and distributed photovoltaic generation, coupled with limited battery adoption, where demand flexibility is sensitive to weather conditions but intertemporal smoothing through storage is constrained.
    \item High participation rates of smart home energy management systems (HEMS), which enable responsive and price-sensitive load adjustments in reaction to differentiated price signals.
    \item Extended training horizons spanning multiple years, particularly when initialized with warm-start parameters differentiated by day type, allowing the algorithm to refine cluster-specific pricing policies over time.
\end{enumerate}
In such settings, incorporating contextual information into the learning process allows the pricing rule to better track context-dependent optimal signals, thereby improving risk reduction relative to homogeneous approaches.

On the other hand, the proposed contextual framework involves additional parameter tuning and may require a longer convergence period. As in most iterative learning algorithms, the selection of the step size remains a central challenge, as it directly affects stability and convergence speed. In our experiments, we find that normalizing the daily demand profile in gradient step substantially improves numerical stability and mitigates sensitivity to the choice of step size. While the remaining hyperparameters are specific to the soft-clustering component, they can be calibrated offline using historical weather data in conjunction with a simulation environment, thereby reducing the burden of online tuning.

In our simulation studies, the cluster-based pricing algorithm typically converges within approximately 10–15 sampled days. However, its performance may exhibit moderate short-term volatility when confronted with previously unobserved day types that differ substantially from the historical regimes used in training. By contrast, the homogeneous feedback-based algorithm generally converges within roughly 10 sampled days to a stable pricing policy. The reasoning behind is that homogeneous feedback-based algorithm exhibits higher price stickiness, as the same price vector is updated uniformly across all days. This reduces day-to-day variability in the price signal and accelerates apparent convergence, even though the resulting policy may be suboptimal under heterogeneous day types. By contrast, the cluster-based contextual algorithm adjusts prices according to cluster probabilities, causing the price signal to vary more dynamically in response to contextual forecasts and consequently requiring a longer period to stabilize.

\begin{figure}[ht]
    \centering
    \includegraphics[width=1\linewidth]{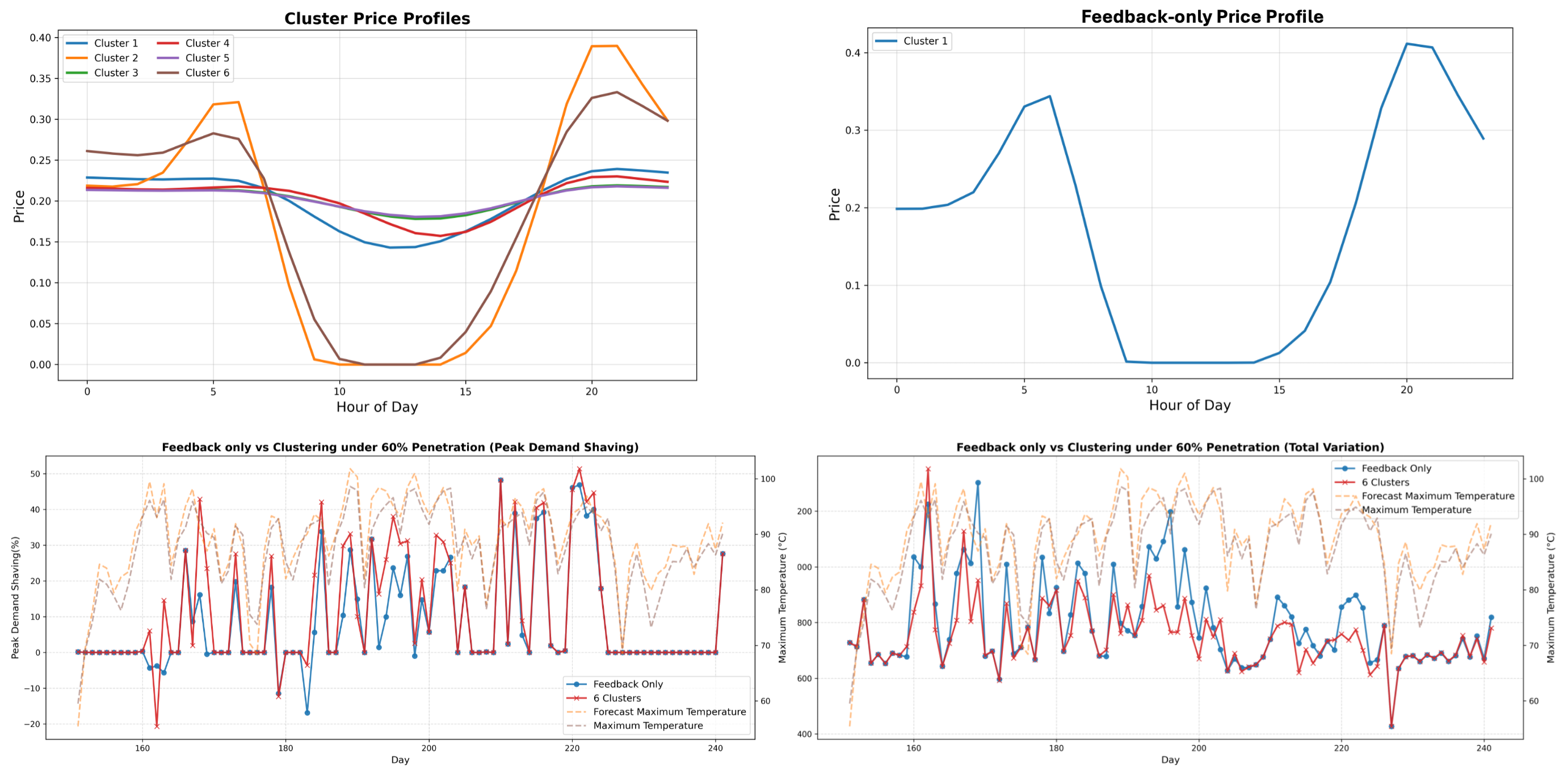}
    \caption{Price Profiles, Peak Demand Shaving and Total Variation\\ under Soft-Clustering and Feedback-Only Algorithms}
    \label{fig:extreme_case}
\end{figure}

Finally, we present a comparative demonstration between (a) the cluster-based soft-clustering algorithm and (b) the homogeneous feedback-only algorithm. The comparison is conducted under a scenario characterized by extreme day-to-day demand variability, including frozen flexible loads and battery states, allowance for power backflow into the grid with symmetric buy and sell prices, full participation of household energy management systems, and 60\% renewable penetration. \autoref{fig:extreme_case} highlights the advantages of the cluster-based approach under these conditions. By explicitly leveraging the cluster structure, the algorithm generates a more diverse set of price profiles, enabling greater peak demand reduction and lower daily demand variability compared to the homogeneous method. Based on these observed benefits, we extend the analysis to multiple case studies and present the corresponding numerical results using the cluster-feedback-based algorithm in Section \ref{sec:numerical}.

\section{Numerical Results}\label{sec:numerical}

To empirically validate the theoretical framework established in Section II and the signal computation algorithms outlined in Section III, this section presents a comprehensive numerical study. Moving beyond stylized analytical models, we implement the proposed one-way dynamic signaling mechanism  following the flowchart in Fig.~\ref{fig:flowchart} within a high-fidelity simulation environment. We deploy the framework on a realistic distribution network populated with hundreds of heterogeneous HEMS, each optimizing a diverse portfolio of behind-the-meter assets—including HVAC units, flexible loads, battery storage, and rooftop photovoltaics—against exogenous weather data and the utility's broadcasted prices. The ensuing subsections detail the test system setup and define key evaluation metrics. We then evaluate the framework's nominal performance, specifically its ability to mitigate device synchronization, achieve Peak Demand Shaving (PDS), and smooth overall load volatility. Finally, we conduct extensive sensitivity analyses to characterize the system's robustness and scalability across varying participation rates, consumer elasticities, geographic climates, and renewable penetration levels.

\subsection{Test System Setups}

\subsubsection{Distribution System} 
To ensure the physical fidelity of our simulations, we utilize an 84-bus (195-node) unbalanced three-phase radial distribution network sourced from the SMART-DS (Synthetic Models for Advanced, Realistic Testing: Distribution Systems and Scenarios) dataset \cite{nrel2023smartds}. This specific network topology was selected because it is highly representative of a typical North American suburban distribution feeder, capturing realistic spatial distributions of line impedances, transformer ratings, and voltage drops. We allocate a total of 486 residential households, each equipped with a HEMS, across the secondary nodes of this network. The households are distributed to mimic natural suburban load densities, ensuring that localized network constraints and aggregate substation behaviors are accurately reflected during the demand response optimization.

\subsubsection{Utility Objective Functions}
The utility cost function is configured following Eq.~\eqref{eq:costfunction} as:
\begin{equation}
    0.9 \cdot QV(\bm{D}) + 0.1 \cdot \|\bm{D}\|^2_2, \nonumber
\end{equation}
which prioritizes smoothing the total demand profile rather than strictly minimizing absolute demand levels. By heavily weighting the variance term, this formulation directly aligns with critical utility objectives: mitigating steep ramping requirements, avoiding secondary ``timer peaks'' caused by synchronized device responses, and reducing reliance on expensive, high-emission peaking generators. Furthermore, while $QV(\bm{D})$ alone is not strongly convex, the inclusion of the strongly convex $L_2$-norm term, $\|\bm{D}\|^2_2$, functions as a crucial mathematical regularization. This addition not only penalizes excessively high overall demand but also significantly enhances the numerical stability and convergence properties of the optimization algorithm.

\subsubsection{HEMS and Devices Setup} 
To ensure realistic household behavior and physical fidelity, the operational parameters and boundary conditions for all controllable assets are derived from industry-standard appliance specifications and typical residential consumption patterns. 
All HEMS are equipped with one HVAC system and one flexible load. To model distributed energy resource integration, we assign a low penetration level of 20\% of HEMS with a combination of rooftop PV and battery. To introduce realistic heterogeneity across the network, individual household parameters are randomized around the following \textbf{average} settings, which directly correspond to the physical models formulated in Section~\ref{sec:devicemodel}:

\textbf{HVAC}: The maximum power rating is bounded at $\overline{p}_{i,d} = 3$~kW. The user's preferred temperature is $T_{i,d}^{prefer} = 75^\circ$F, tightly constrained by upper and lower comfort bounds of $[\underline{T}_{i,d}, \overline{T}_{i,d}] = [72^\circ\text{F}, 78^\circ\text{F}]$. The thermal dynamic parameters, specifying cooling effectiveness and house insulation properties, are averaged at $\zeta_1 = 0.9$ and $\zeta_2 = 1.0$, respectively.

\textbf{Flexible Load}: The preferred hourly load profile, $p_{i,d}^{prefer}$, is predefined based on baseline residential data. The operational flexibility bounds are set to $\pm 20\%$ of these preferred hourly values, with peak-hour flexibility reduced to about $\pm 10\%$. We strictly enforce that the total daily energy consumption $E_{i,d}$ for the flexible load remains fixed over the 24-hour horizon, consistent with the energy conservation constraints of the model.

\textbf{PV}: The solar inverter is sized with a maximum power rating of 5~kW. Following the generation sign convention established in our framework, the output power is constrained within $\overline{p}_{i,d} \le p_{i,d} \le 0$~kW, where $\overline{p}_{i,d} = -5$~kW.

\textbf{Battery}: The maximum charging and discharging power ratings are constrained by $\overline{p}_{i,d} = 5$~kW and $\underline{p}_{i,d} = -5$~kW, respectively. The total energy capacity is set such that the battery can continuously charge or discharge at its maximum power rating for approximately 4 hours before reaching its operational limits. Furthermore, the preferred state-of-charge is defined as $SOC_{i,d}^{prefer} = 50\%$, and the operational SOC is securely bounded within $[\underline{SOC}_{i,d}, \overline{SOC}_{i,d}] = [20\%, 80\%]$ to preserve battery health and longevity.

\textbf{No-sell Constraints}: We impose a strict non-export limit, assuming that households do not inject reverse power flow back into the distribution grid. Consequently, the aggregate net power consumption for any household $i$ must remain non-negative at all times:
\begin{equation}
    p_{i}=\sum_{d\in\mathcal{D}_i} p_{i,d} \ge 0. \nonumber
\end{equation}

Given these operational setup parameters and physical constraints, the local HEMS co-optimization problem—defined by \eqref{eq:opt1} for participating households and \eqref{eq:opt2} for non-participating households—is solved using the interior-point algorithm provided by the IPOPT solver.

\subsubsection{Temperature and Irradiance}\label{sec:cityinputs}

\begin{figure}[!h]
\centerline{\includegraphics[width=.3\textwidth]{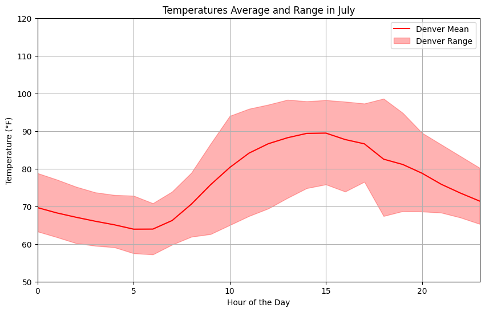}
\includegraphics[width=.3\textwidth]{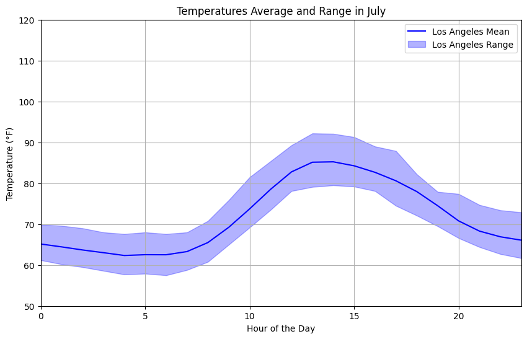}
\includegraphics[width=.3\textwidth]{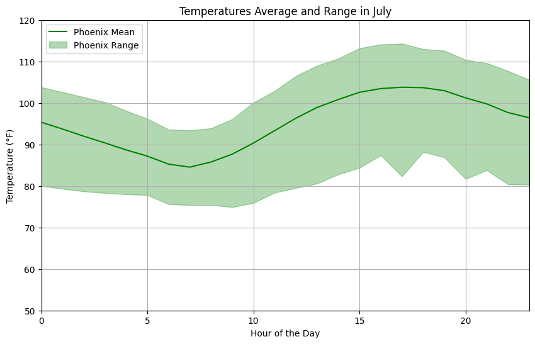}}
\centerline{\includegraphics[width=.3\textwidth]{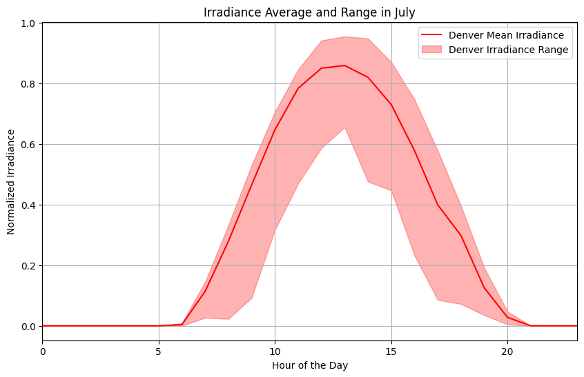}
\includegraphics[width=.3\textwidth]{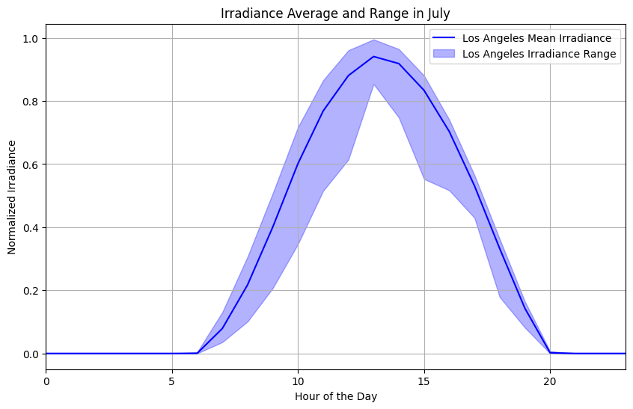}
\includegraphics[width=.3\textwidth]{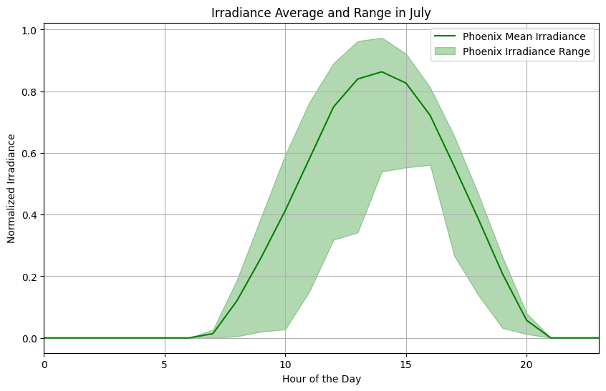}}
\caption{(Upper) Daily average temperatures and temperature ranges for Denver, LA, and Phoenix; (lower) daily average solar irradiance and irradiance ranges for these three cities.}
\label{fig:cities_inputs}
\end{figure}

We collect forecasted and realized historical daily temperature and solar irradiance data from Open-Meteo \cite{zippenfenig2023open}. Specifically, we query the hourly variables $\mathtt{temperature\_2m}$ and $\mathtt{Shortwave\_Solar\_Radiation\_GHI (Instant)}$. To accurately reflect the realized weather conditions, we utilize the historical hourly records stored in the ECMWF Integrated Forecasting System (IFS). 
We collect and use the 2022 Denver area temperature and solar irradiance data for all case studies, except in Section~\ref{sec:cities} where we explore different cities and weather patterns. 
The summer temperature and solar irradiance profiles for the three cities under study are illustrated in Fig.~\ref{fig:cities_inputs}.

\subsubsection{Benchmark Scenario}
To establish a baseline for evaluating the proposed dynamic pricing framework, we define a price-oblivious benchmark scenario where the participation rate is set to zero ($\mathcal{P}=0$). In this uncoordinated state, all 486 HEMS operate independently, optimizing energy consumption solely to satisfy user preferences and internal device constraints, without any regard for external utility signals or grid conditions. Mathematically, this implies that the entire population belongs to the non-participating set ($\mathcal{N}_{out}$), with each HEMS individually determining its operational profile by solving the local optimization problem defined in \eqref{eq:opt2}.

\subsubsection{Nominal Case Setup}
To evaluate the framework under a realistic, partial-adoption setting, we define a nominal scenario with a participation rate of $\mathcal{P}=2/3$. Under this configuration, 322 HEMS opt into the demand flexibility program, forming the participating subset $\mathcal{N}_{in}$. These active households determine their device schedules by solving the price-responsive optimization problem defined in \eqref{eq:opt1}. The remaining 164 HEMS constitute the non-participating subset $\mathcal{N}_{out}$; they remain price-oblivious and dictate their device outputs strictly based on user preferences by solving \eqref{eq:opt2}. This heterogeneous mix of active and passive households accurately reflects the practical realities of consumer adoption within a distribution network.

\subsubsection{DF Performance Limits: Two-Way Communication and Direct Control} \label{sec:defineidealcases}
To rigorously evaluate the efficacy and optimality gap of the proposed one-way signaling framework, we define two idealized reference scenarios. While both scenarios are challenging to implement in practice due to infrastructure and privacy constraints, they provide crucial theoretical upper bounds for system performance:

\begin{itemize}
    \item \textbf{Two-Way Communication Case:} In contrast to the proposed one-way framework—where the utility updates signals based solely on the realized demand of the previous 24 hours—a two-way communication architecture allows for iterative, day-ahead negotiations. In this scenario, we allow the utility and the HEMS to repeatedly exchange provisional price signals and anticipated demand profiles. By conducting sufficient iterations until the system reaches convergence, this approach generates a theoretically optimal demand response. The converged results serve as a rigorous performance bound to quantify the efficiency of our one-way communication algorithm.

    \item \textbf{Direct Load Control (DLC) Case:} Rather than fundamentally restructuring the mathematical model to accommodate centralized utility scheduling, we simulate an ideal DLC paradigm by modifying the behavioral parameters within our existing framework. Specifically, we scale the user preference elasticity parameters, $\gamma_{i,d}$, to approach zero for all controllable devices. Under this condition, the HEMS mathematically disregards user comfort (though strictly adhering to physical device limits) and optimizes solely for utility signal responsiveness. When this parameter adjustment is combined with the iterative two-way communication process described above, the setup effectively mimics a centralized direct control scenario, yielding the absolute theoretical maximum for grid-level flexibility.
\end{itemize}

\subsection{Evaluation Metrics}
To systematically evaluate the efficacy of the proposed demand flexibility framework, we define the following key performance metrics. Let $\bm{D}_0$ denote the aggregated demand vector of the uncoordinated benchmark scenario, and $\bm{D}$ denote the realized demand vector of the control scenario under evaluation.

\subsubsection{Peak Demand Shaving Percentage (PDS\%)}
We define the daily PDS\% as the percentage reduction in the daily maximum load compared to the benchmark:
\begin{equation}
    \text{PDS}\%(\bm{D}|\bm{D}_0) = \frac{\max(\bm{D}_0) - \max(\bm{D})}{\max(\bm{D}_0)} \times 100\%.
\end{equation}
A positive PDS\% indicates a successful reduction in peak demand, whereas a negative value implies an unintended exacerbation of the peak load.

\subsubsection{Monthly Peak Demand Shaving (MPS\%) \& Aggregated Monthly Peak Demand Shaving (AMPS\%)}
To evaluate performance over an extended period, we define two monthly metrics. MPS\% evaluates the absolute maximum peak reduction over an entire calendar month. It is computed by comparing the single highest peak of the evaluation demand profile across the month against the single highest peak demand of the benchmark scenario. Given the daily evaluation demand $\bm{D}(d)$ and benchmark $\bm{D}_0(d)$ for day $d \in \{1, \ldots, N_m\}$ (where $N_m$ is the number of days in the month), MPS\% is defined as:
\begin{equation}
    \text{MPS}\% = \frac{\max_{d} \big(\max(\bm{D}_0(d))\big) - \max_{d} \big(\max(\bm{D}(d))\big)}{\max_{d} \big(\max(\bm{D}_0(d))\big)} \times 100\%. \nonumber
\end{equation}
While MPS\% evaluates the absolute worst-case peak of the month, the AMPS\% evaluates the aggregate daily peak reduction. AMPS\% is computed as the ratio of the total shaved peak load to the total baseline peak load across all days in the month:
\begin{equation}
    \text{AMPS}\% = \frac{\sum_{d=1}^{N_m} \Big(\max(\bm{D}_0(d)) - \max(\bm{D}(d))\Big)}{\sum_{d=1}^{N_m}\max(\bm{D}_0(d))} \times 100\%. \nonumber
\end{equation}

\subsubsection{Maximum Hourly Variation (Ramp Rate)}
To evaluate the smoothness of the demand profile, we measure the maximum absolute change in demand from one hour to the next. To distinguish this metric from the temporal variance term used in the utility objective function, we denote this maximum hourly variation as $\Delta_{\max}(\bm{D})$:
\begin{equation}
    \Delta_{\max}(\bm{D}) = \max_{t=1, \ldots, T-1}(|D_{t+1}-D_t|).
\end{equation}
This metric characterizes the steepest ramping requirement imposed on the grid. A smaller $\Delta_{\max}$ indicates a smoother demand profile, with a perfectly flat profile yielding a value of zero.

\subsubsection{Load Factor (LF)}
The Load Factor (LF) falling between 0 and 1 assesses the utilization efficiency of the grid infrastructure over the 24-hour horizon $T$, defined as:
\begin{equation}
    \text{LF}(\bm{D}) = \frac{\sum_{t=1}^T D_t}{\max_t(D_t)\times T}. \nonumber
\end{equation}
The numerator represents the total energy consumed, while the denominator represents the theoretical maximum energy consumption if the peak demand were sustained continuously. A higher LF indicates a flatter demand profile, whereas a lower LF characterizes a highly volatile or ``spiky'' load.

\subsubsection{Total Energy Reduction}
To confirm that the primary mechanism of our demand response framework is load \textit{shifting} rather than load \textit{shedding}, we track the total daily energy reduction. A healthy flexibility response should maintain consumer utility by keeping total daily energy consumption relatively constant while effectively displacing it away from peak hours.

\subsection{Nominal Case Results}

\subsubsection{Result for A Typical Day} 
\begin{figure}[!h]
\centerline{\includegraphics[width=.6\textwidth]{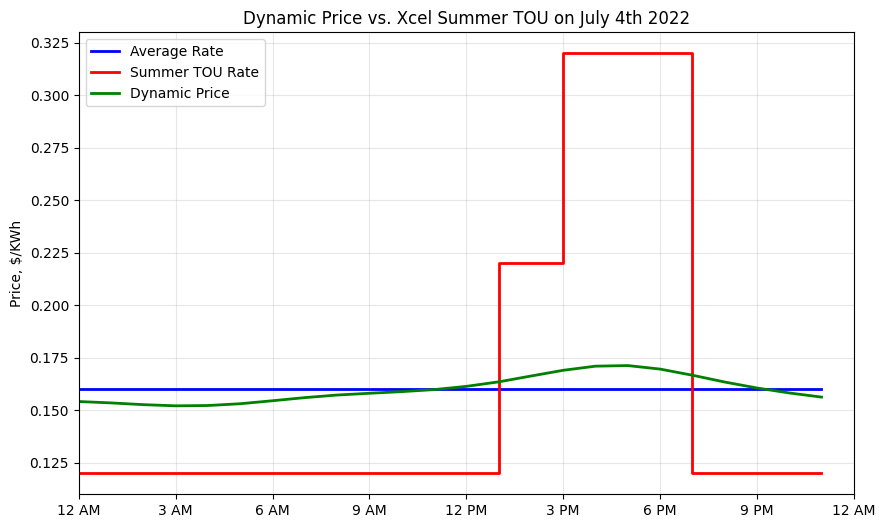}}
\caption{Prices on July 4th, 2022. {\color{red}Add a zoomed-in figure.}}
\label{fig:typicalprices}
\end{figure}

\begin{figure}[!h]
\centerline{\includegraphics[width=.4\textwidth]{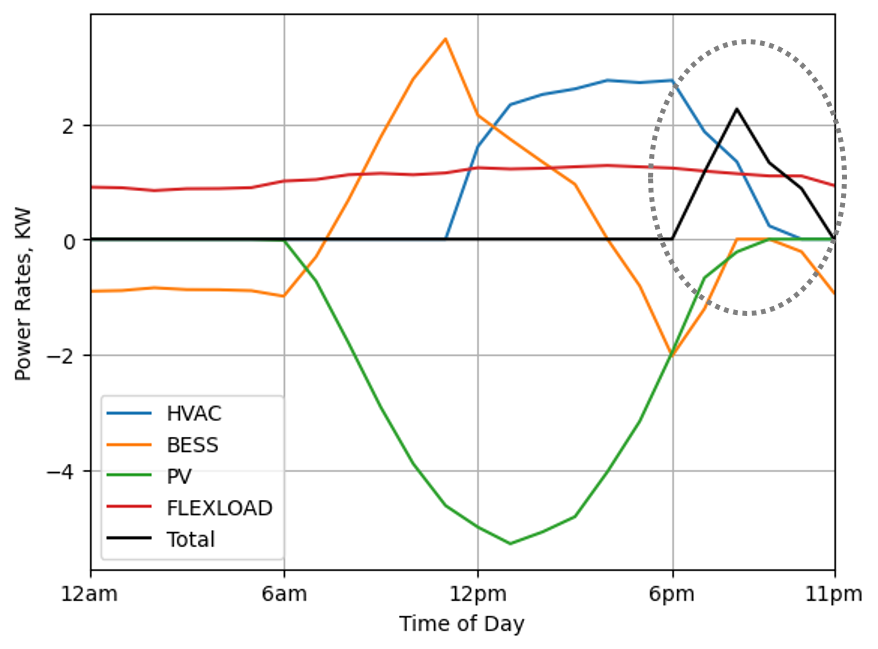}
\includegraphics[width=.4\textwidth]{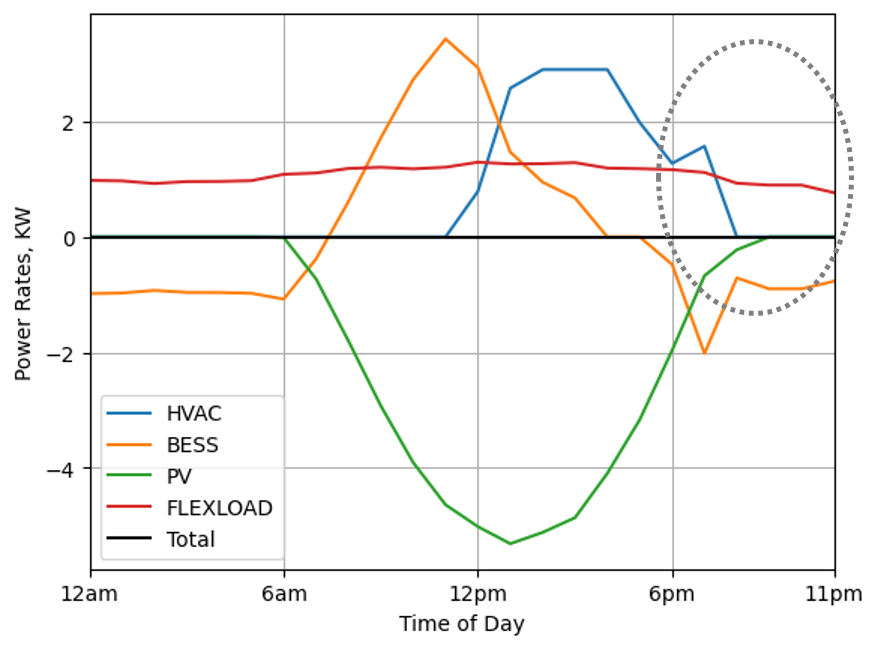}}
\caption{Daily demand profiles for a household from benchmark (left) and under dynamic price (right) on July 4th, 2022. The circled hours show improvement of the total household power usage thanks to dynamic pricing.}
\label{fig:hemsdailydemands}
\end{figure}

\begin{figure}[!h]
\centerline{\includegraphics[width=.7\textwidth]{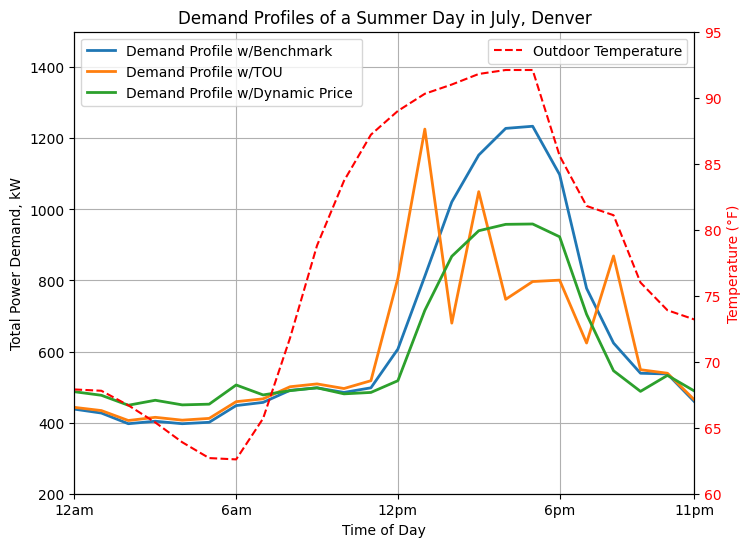}}
\caption{Daily demand profiles from benchmark and under dynamic price on July 4th, 2022.}
\label{fig:dailydemands}
\end{figure}

To illustrate the immediate impact of the proposed mechanism, we first examine the load profile of a typical summer day. Fig~\ref{fig:typicalprices} shows the prices on this day and Fig.~\ref{fig:hemsdailydemands} compares the demand of a single HEMS under benchmark case and under dynamic price. The results demonstrate that the HEMS successfully shifts discretionary loads—primarily HVAC setpoints—away from peak price periods: the output from this particular household becomes net zero. Such individual response aggregates into a significantly flattened system-wide demand curve as shown in Fig.~\ref{fig:dailydemands}, providing a baseline for the year-wide statistical analysis that follows. Here, the demand with TOU are based on summer rates from Fig~\ref{fig:xceltou}.

\subsubsection{PDS\%, MPS\%, and Load Variation}
\begin{figure}[!h]
\centerline{\includegraphics[width=.7\textwidth]{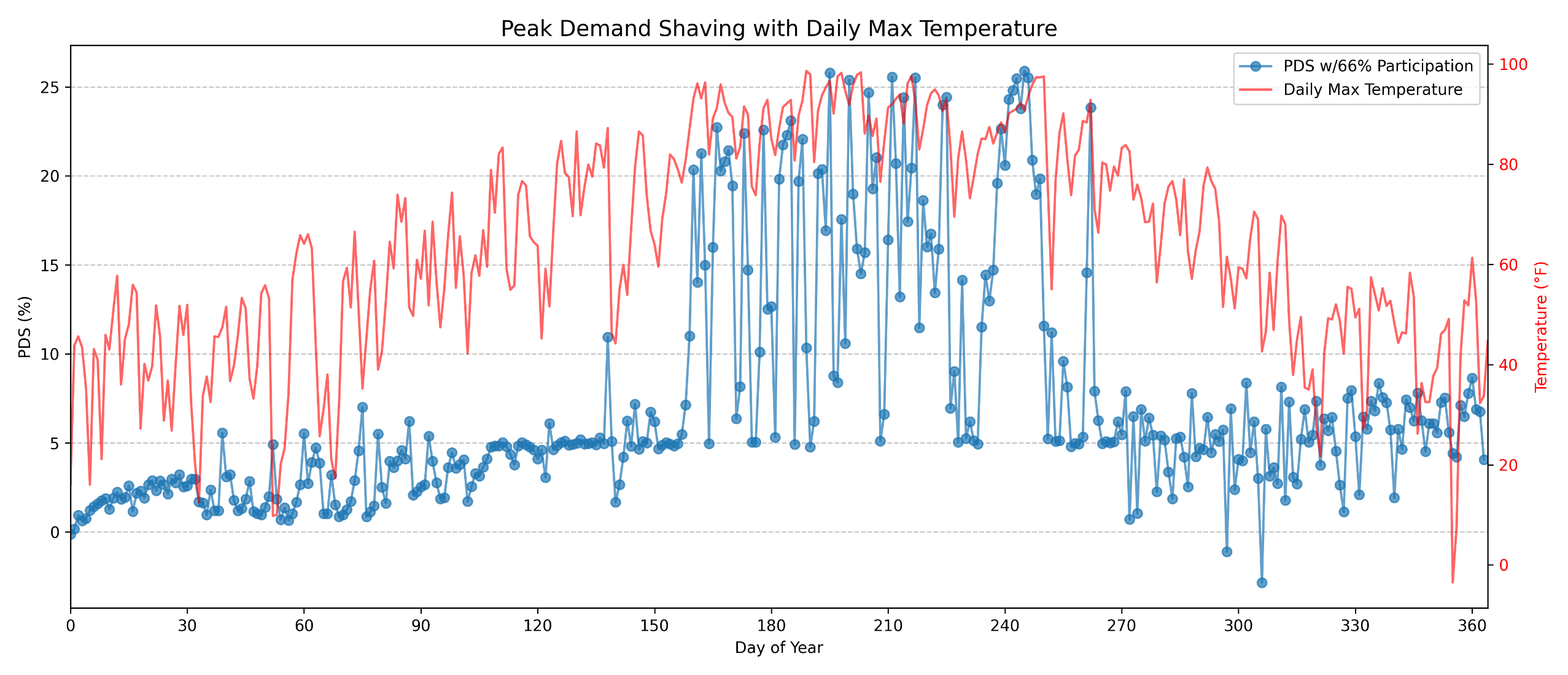}}
\centerline{\includegraphics[width=.7\textwidth]{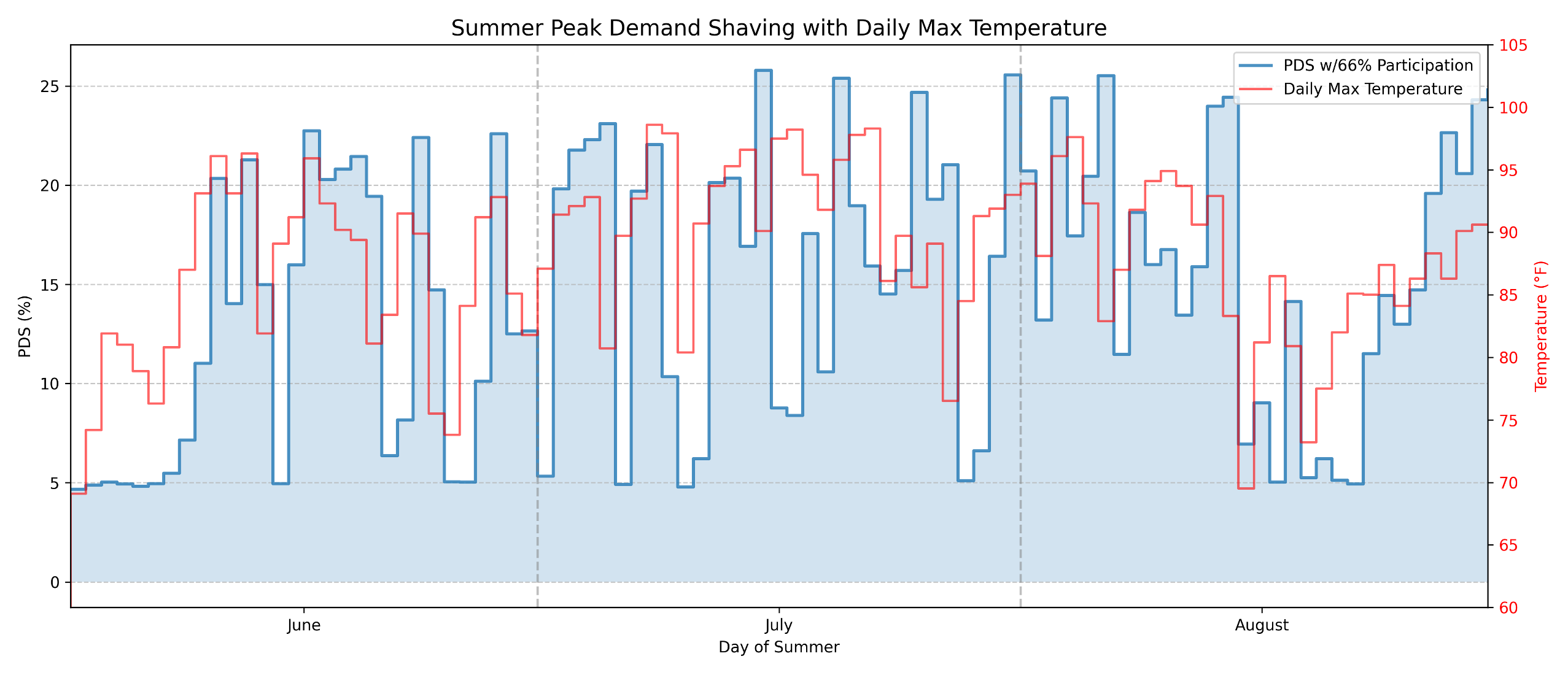}}
\caption{Daily PDS improvement for a year (upper), and for summer months (lower).}
\label{fig:nominal_PDS}
\end{figure}

Fig.~\ref{fig:nominal_PDS} (upper) illustrates the daily Peak Demand Shaving (PDS\%) of the nominal case relative to the benchmark over a full calendar year, overlaid with maximum daily temperatures. In Fig.~\ref{fig:nominal_PDS} (lower), we focus on the summer months (June, July, and August), where demand typically peaks due to intensive HVAC usage. 

The results indicate that dynamic pricing reduces daily peak demand consistently throughout the year, with a positive PDS\% achieved on 99\% of days. Notably, during high-temperature summer days, PDS\% exceeds 20\%, compared to approximately 5\% during colder periods when HVAC activity is minimal.

\begin{figure}[!h]
\centerline{
\includegraphics[width=.48\textwidth]{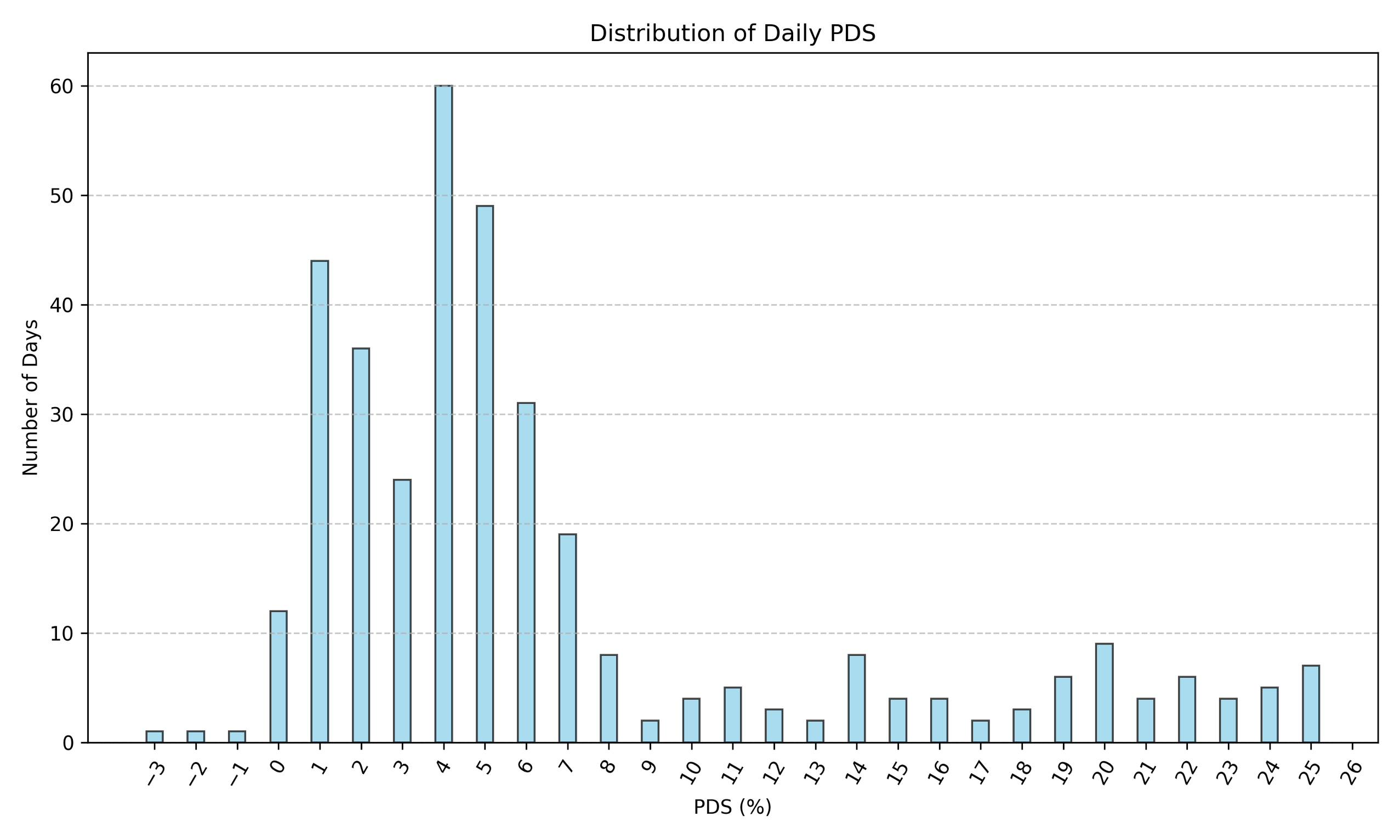}
\includegraphics[width=.49\textwidth]{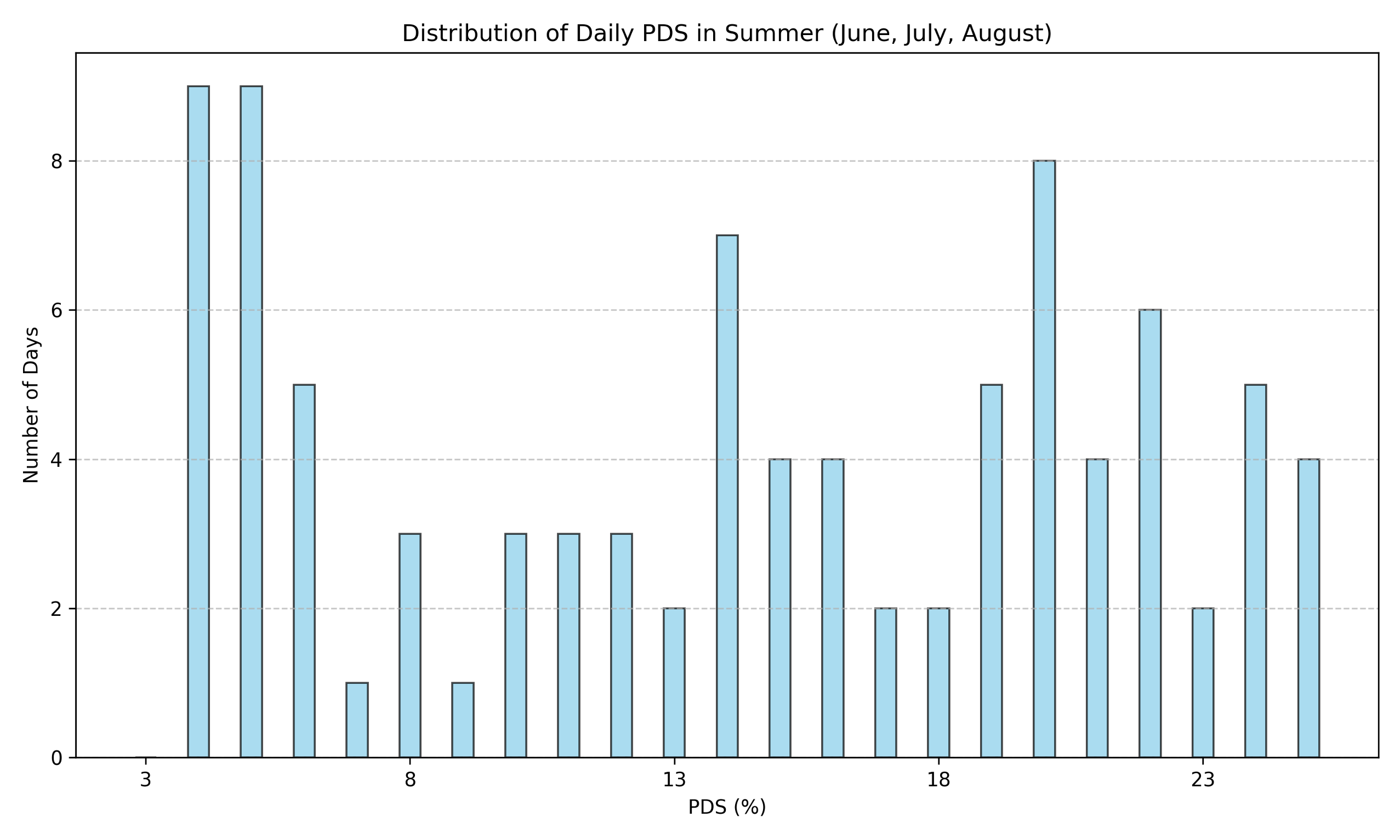}}
\caption{PDS distribution for a year and for summer.}
\label{fig:nominal_PDS_distribution}
\end{figure}

The statistical distribution of these improvements is shown in Fig.~\ref{fig:nominal_PDS_distribution}, while Table~\ref{tab:monthly_mps} provides the Monthly Peak Shaving (MPS\%) values. These metrics clarify how the demand response effectiveness scales with seasonal load intensity. Furthermore, the reduction in total demand variation (load volatility) is calculated at 11.12\% for the full year and 19.33\% specifically for the summer, reflecting a more stable and predictable grid profile.

\begin{table}[h!]
    \centering
    \caption{MPS\% and AMPS\% Data.}
    \label{tab:monthly_mps}
    \begin{threeparttable}
    \begin{tabular}{lcccccc}
        \toprule
        \textbf{Month} & \textbf{Jan} & \textbf{Feb} & \textbf{Mar} & \textbf{Apr} & \textbf{May} & \textbf{Jun} \\
        \midrule
        \textbf{MPS\%} & 0 & 2.96 & 4.73 & 4.44 & 6.72 & 14.02 \\
        \midrule
        \textbf{AMPS\%} & 1.90 & 1.96 & 3.05 & 3.85 & 5.10 & 14.33 \\
        \midrule
        \textbf{Month} & \textbf{Jul} & \textbf{Aug} & \textbf{Sep} & \textbf{Oct} & \textbf{Nov} & \textbf{Dec} \\
        \midrule
        \textbf{MPS\%} & 10.33 & 15.16 & 18.90 & 2.38 & 2.54 & 7.11 \\
        \midrule
        \textbf{AMPS\%} & 16.67 & 16.89 & 13.54 & 4.65 & 4.68 & 6.33 \\
        \bottomrule
    \end{tabular}
    \end{threeparttable}
    \end{table}

\subsubsection{Load Factor (LF)} 
\begin{figure}[!h]
\centerline{\includegraphics[width=.7\textwidth]{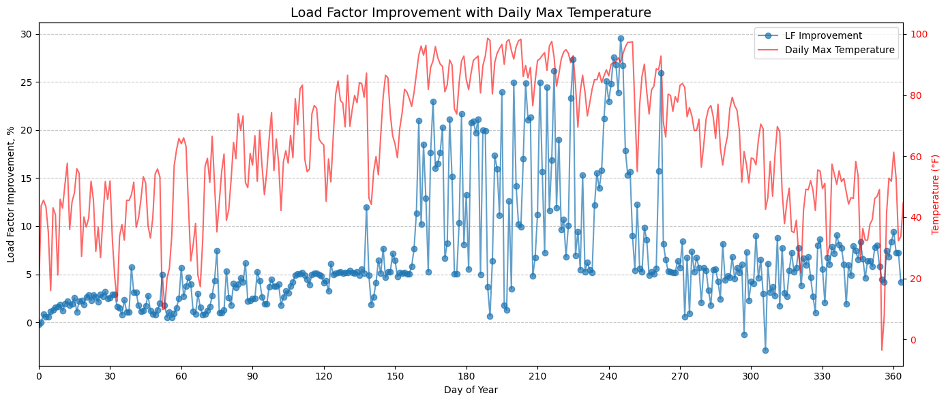}}
\centerline{\includegraphics[width=.7\textwidth]{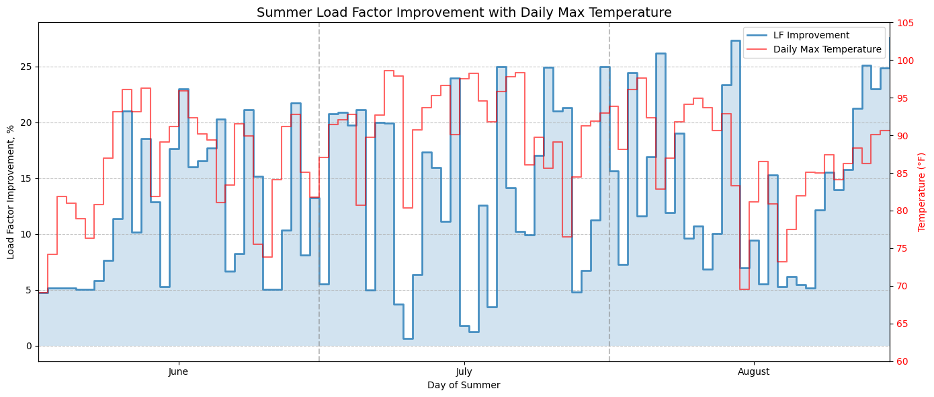}}
\caption{LF Improvement for a year (upper) and over summer months (lower).}
\label{fig:nominal_LF}
\end{figure}

\begin{figure}[!h]
\centerline{\includegraphics[width=.4\textwidth]{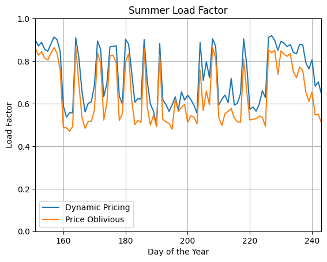}}
\caption{LF of nominal case and benchmark.}
\label{fig:nominal_LF_summer}
\end{figure}

The daily Load Factor (LF) improvement is presented in Fig.~\ref{fig:nominal_LF}, with a specific focus on summer performance in Fig.~\ref{fig:nominal_LF_summer}. We observe that summer days with an initial LF of 0.5 are improved to over 0.6—a relative increase of more than 20\%. This improvement signifies a less "spiky" demand profile, resulting from the effective shifting of flexible demand away from peak hours.

\subsubsection{Energy Usage vs. Peak Reduction}
The impact of the nominal case on total energy consumption is summarized below:
\begin{itemize}
    \item \textbf{Annual:} 1.4\% reduction in total energy usage.
    \item \textbf{Summer:} 4.2\% reduction in total energy usage.
\end{itemize}
These values are significantly lower than the PDS improvements (9\% annually and 16\% in summer). This discrepancy confirms that the primary driver of peak reduction is \textit{load shifting} rather than \textit{load shedding}. This is a critical finding, as it suggests the mechanism achieves grid stability without requiring substantial customer discomfort or causing significant revenue loss for the utility.

\subsubsection{Comparison with Two-Way Communication and Direct Control}
To evaluate the efficiency of the proposed system, the nominal case is compared against the ideal benchmarks defined in Section~\ref{sec:defineidealcases}. While full visual results for these cases are provided in the Appendix due to their similarity to the nominal case, the comparative summer performance is summarized as follows:
\begin{itemize}
    \item \textbf{Two-way communication} achieves 17.8\% PDS.
    \item \textbf{Direct control} (upper bound) achieves 19.0\% PDS.
\end{itemize}
These results demonstrate that our proposed one-way communication-based mechanism, despite its lower implementation complexity and higher privacy, achieves performance remarkably close to ideal, centralized control scenarios.

\subsubsection{Section Summary} 
The analysis of the nominal case demonstrates that a feedback-based, day-ahead dynamic pricing structure effectively incentivizes HEMS-equipped households to optimize their energy patterns. This results in a smoothed demand profile and reduced peak levels without compromising consumer comfort. Furthermore, the performance gap between this practical one-way structure and ideal control methods is minimal, validating the feasibility of the proposed approach.

\subsection{Study of Important System Setup Parameters}
In this part, we fixed all system parameters in the nominal cases but one under study for each subsection. These results help us understand how different parameters affect the demand response results, from system's perspective.

\subsubsection{Participation Rates}

\begin{figure}[!h]
\centerline{\includegraphics[width=.3\textwidth]{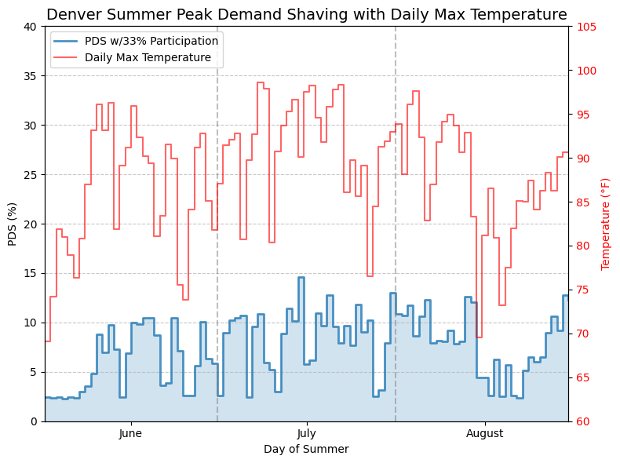}
\includegraphics[width=.3\textwidth]{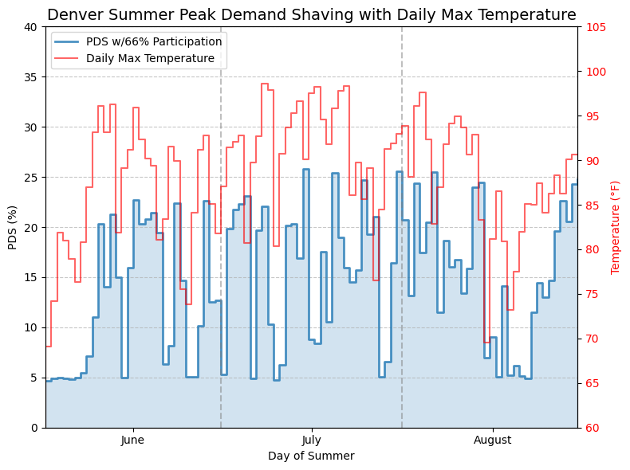}
\includegraphics[width=.3\textwidth]{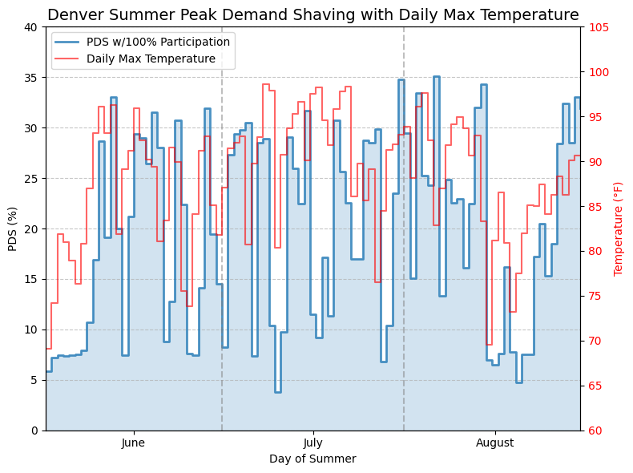}
}
\caption{Summer PDS improvement under different participation rates. From left to right: 33\%, 66\%, and 100\% participation.}
\label{fig:participation}
\end{figure}

To evaluate the scalability of the proposed mechanism, we vary the participation rate from 33\%, to 66\% (the nominal case), and finally to 100\%. The resulting summer PDS\% improvements are illustrated in Fig.~\ref{fig:participation}, with a detailed summary of PDS\% and variation reduction provided in Table~\ref{tab:participation}. 

As anticipated, increasing the participation rate from 0\% (the benchmark) to 66\% yields a near-linear improvement in PDS\%. However, as participation approaches 100\%, the rate of improvement begins to exhibit diminishing marginal returns. A similar trend is observed in the reduction of demand variation. This phenomenon occurs because further flattening an already smoothed demand profile becomes increasingly difficult; in an extreme example, no further peak reduction is possible once the demand profile reaches a perfectly flat line.

Therefore, from a practical implementation standpoint, these results suggest that there is no need to pursue an extremely high participation rate. Since the majority of the grid-balancing benefits are captured at moderate participation levels (e.g., 66\%), the marginal utility of enrolling the remaining 33\% of customers may not justify the additional marketing and infrastructure costs required to achieve universal adoption.

\begin{table}[h]
    \centering
    \caption{Participation Rates, PDS\%, and Variation Reduction for summer months.}\label{tab:participation}
    \begin{tabular}{lll}
        \toprule
        \textbf{Participation levels} & \textbf{PDS\%} & \textbf{Variation Reduction \%} \\
        \midrule
        33\% & 8.27  & 12.17 \\
        66\% & 16.01 & 19.33 \\
        100\%  & 21.42 & 19.60 \\
        \bottomrule
    \end{tabular}
\end{table}

\subsubsection{Sensitivity Analysis of Elasticity Coefficients}

\begin{figure}[!h]
\centerline{\includegraphics[width=.4\textwidth]{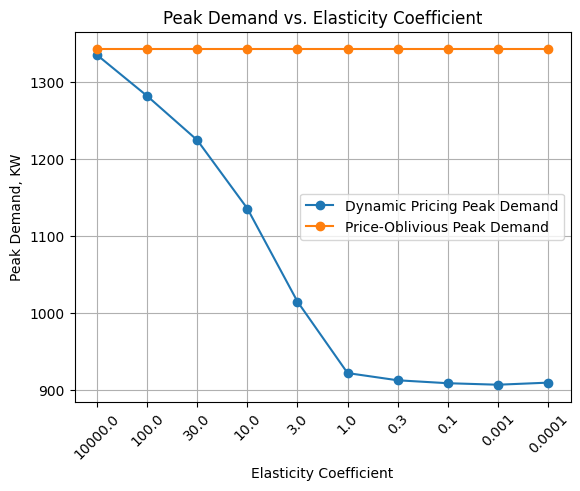}}
\caption{Daily peak demand with respect the elasticity.}
\label{fig:elasticity}
\end{figure}

To evaluate the impact of consumer flexibility on system performance, we conduct a sensitivity analysis by scaling the elasticity parameters $\gamma_{i,d}$ for all participating devices. The scaling factor ranges from 10,000 to 0.0001, effectively shifting the optimization priority between user comfort and price responsiveness. Fig.~\ref{fig:elasticity} illustrates the peak demand for a representative summer day across these varying elasticities.

The results reveal three distinct operational phases:
\begin{itemize}
    \item \textbf{Preference Dominance (Scale $\geq$ 100):} At high scaling factors, the cost function places an overwhelming weight on user preferences. Consequently, HEMS behaves almost identically to a non-participating (price-oblivious) household, with peak demand remaining near the benchmark of 1,340 kW.
    \item \textbf{High Responsiveness (100 $>$ Scale $\geq$ 1):} As the scaling factor decreases, the HEMS becomes significantly more sensitive to dynamic price signals. In this region, we observe a sharp reduction in peak demand as the system aggressively shifts loads to lower-price periods to minimize costs.
    \item \textbf{Constraint Saturation (Scale $<$ 1):} Below a scaling factor of 1, the peak demand reduction reaches an asymptote near 900 kW. Despite the HEMS being mathematically indifferent to user comfort at this stage, further demand lowering is curtailed by two primary factors: (i) the demand profile has already reached a state of maximum feasible smoothness, and (ii) hard operational constraints, such as HVAC temperature deadbands, prevent further load shifting.
\end{itemize}

This analysis demonstrates that there is a practical performance bound for peak demand reduction. Even under extreme price sensitivity, the proposed mechanism respects the physical constraints of the household, ensuring that the demand response remains realistic and never compromises essential device functionality.

\subsubsection{Different Cities/Climates/Weather Patterns}\label{sec:cities}

\begin{figure}[!h]
\centerline{\includegraphics[width=.3\textwidth]{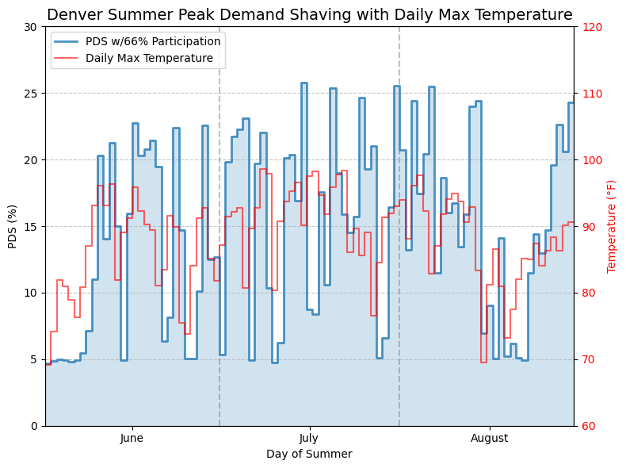}
\includegraphics[width=.3\textwidth]{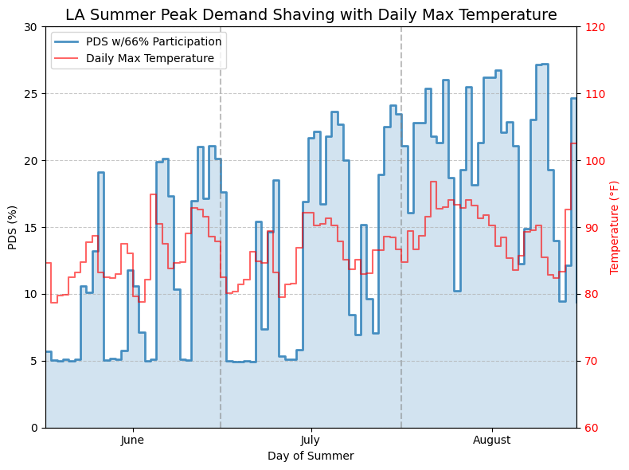}
\includegraphics[width=.3\textwidth]{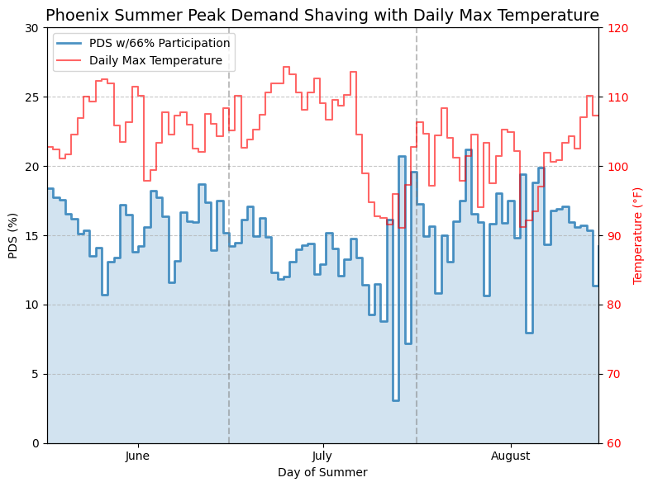}
}
\caption{Summer PDS improvement for three cities.}
\label{fig:cities_PDS}
\end{figure}

To evaluate the robustness of the proposed mechanism across diverse climatic conditions, we vary the solar irradiance and temperature inputs to represent three distinct regions: Denver, Los Angeles, and Phoenix, as detailed in Section~\ref{sec:cityinputs}. The resulting daily summer PDS\% improvements are illustrated in Fig.~\ref{fig:cities_PDS}, with aggregate metrics summarized in Table~\ref{tab:cities}.

The results indicate that the framework remains effective across different climates, though the performance characteristics vary with local weather patterns. Notably, while Phoenix has the highest cooling requirements, it achieves the lowest PDS\% (13.85\%) among the three cities. This is attributed to the "flexibility saturation" that occurs during extreme heat; when ambient temperatures remain persistently high, HVAC systems must operate near continuously to maintain indoor temperatures within the required deadbands, leaving little room for discretionary load shifting.

However, Phoenix demonstrates the most significant improvement in variation reduction (40.66\%). This suggests that while individual peaks are harder to "shave" in extreme climates, the overall demand profile becomes substantially more stable and predictable under the proposed dynamic pricing. In contrast, Denver and Los Angeles show higher PDS\% potential due to more moderate temperature fluctuations, which allow for greater flexibility in scheduling residential loads.

\begin{table}[h]
    \centering
    \caption{City Data: PDS\% and Variation Reduction}\label{tab:cities}
    \begin{tabular}{lll}
        \toprule
        \textbf{City} & \textbf{PDS\%} & \textbf{Variation Reduction \%} \\
        \midrule
        Denver      & 16.01 & 19.33 \\
        Los Angeles & 16.58 & 21.12 \\
        Phoenix     & 13.85 & 40.66 \\
        \bottomrule
    \end{tabular}
\end{table}

\subsubsection{Impact of PV-Battery Penetration Levels}

To investigate the compatibility of the proposed system in a grid with high renewable integration, we increase the penetration level of the PV-battery combination from 20\% (the nominal case) to 60\%. In this analysis, the nominal case—rather than the price-oblivious benchmark—serves as the baseline to isolate the specific impact of increased distributed energy resource (DER) adoption. 

\begin{figure}[!h]
\centerline{\includegraphics[width=.7\textwidth]{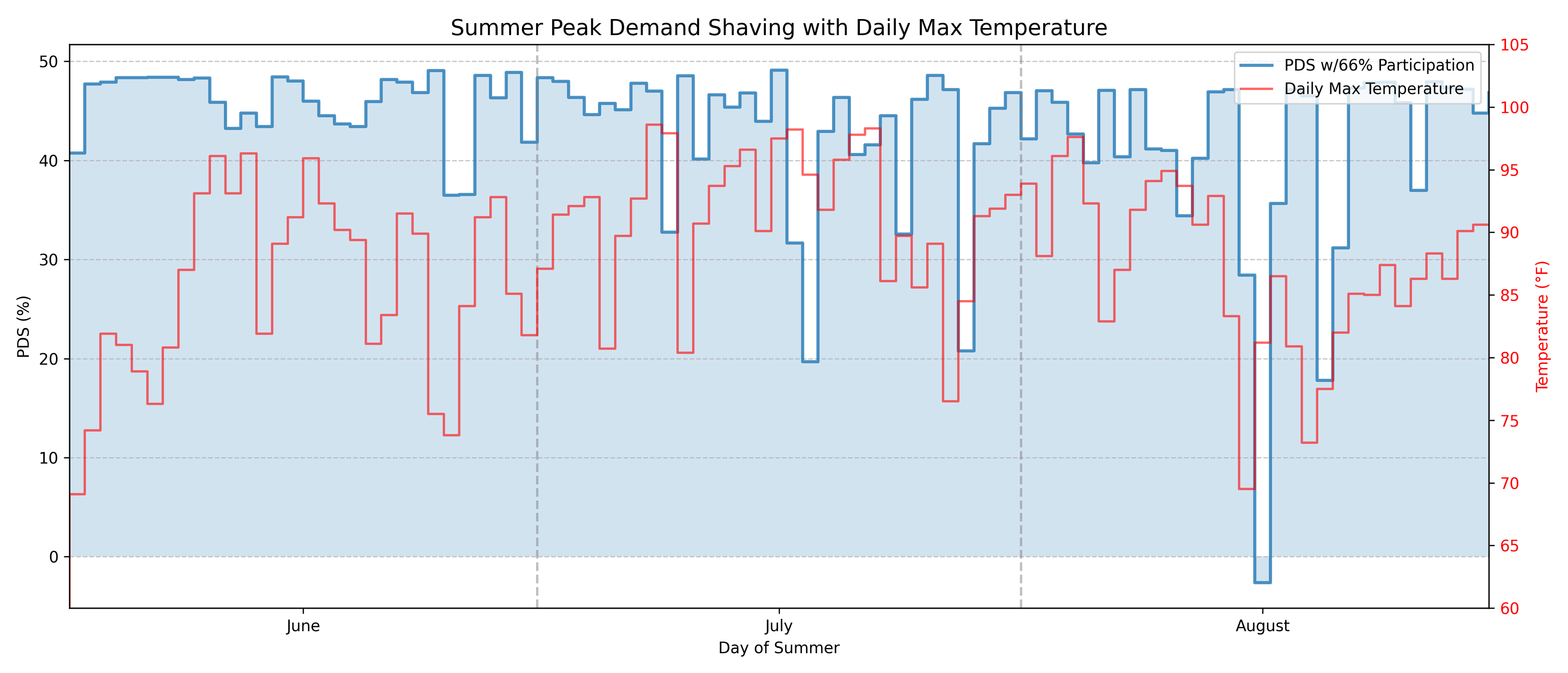}}
\caption{PDS Improvement from 20\% to 60\% penetration of renewable energy over summer months.}
\label{fig:renewable_PDS}
\end{figure}

The summer PDS\% improvement for the 60\% penetration scenario is illustrated in Fig.~\ref{fig:renewable_PDS}. Under this high-penetration regime, the system achieves an overall summer PDS\% of 32.79\% and a variation reduction of 31.00\% relative to the nominal case. This significant performance boost is primarily driven by the temporal alignment between solar generation and summer cooling loads; the additional PV output effectively offsets the original afternoon demand peaks, while the increased battery capacity provides the HEMS with greater flexibility to shift residual loads.

The results confirm that the proposed demand response mechanism is highly complementary to renewable integration. Increasing the penetration of DERs does not lead to diminishing returns in this context; instead, it creates a compounding effect where local generation and storage significantly amplify the peak-shaving and grid-stabilizing capabilities of the HEMS.

\subsubsection{HVAC-Only}

In this subsection, we examine a scenario where only HVAC systems are responsive to utility signals, while all other household appliances remain price-oblivious. This configuration reflects the current industry trend where utilities prioritize using high-consumption cooling and heating units for demand response purpose. To quantify the performance of the HVAC-only case, we compare these results against the nominal case where all HEMS-connected devices participate.

\begin{figure}[!h]
\centerline{\includegraphics[width=.7\textwidth]{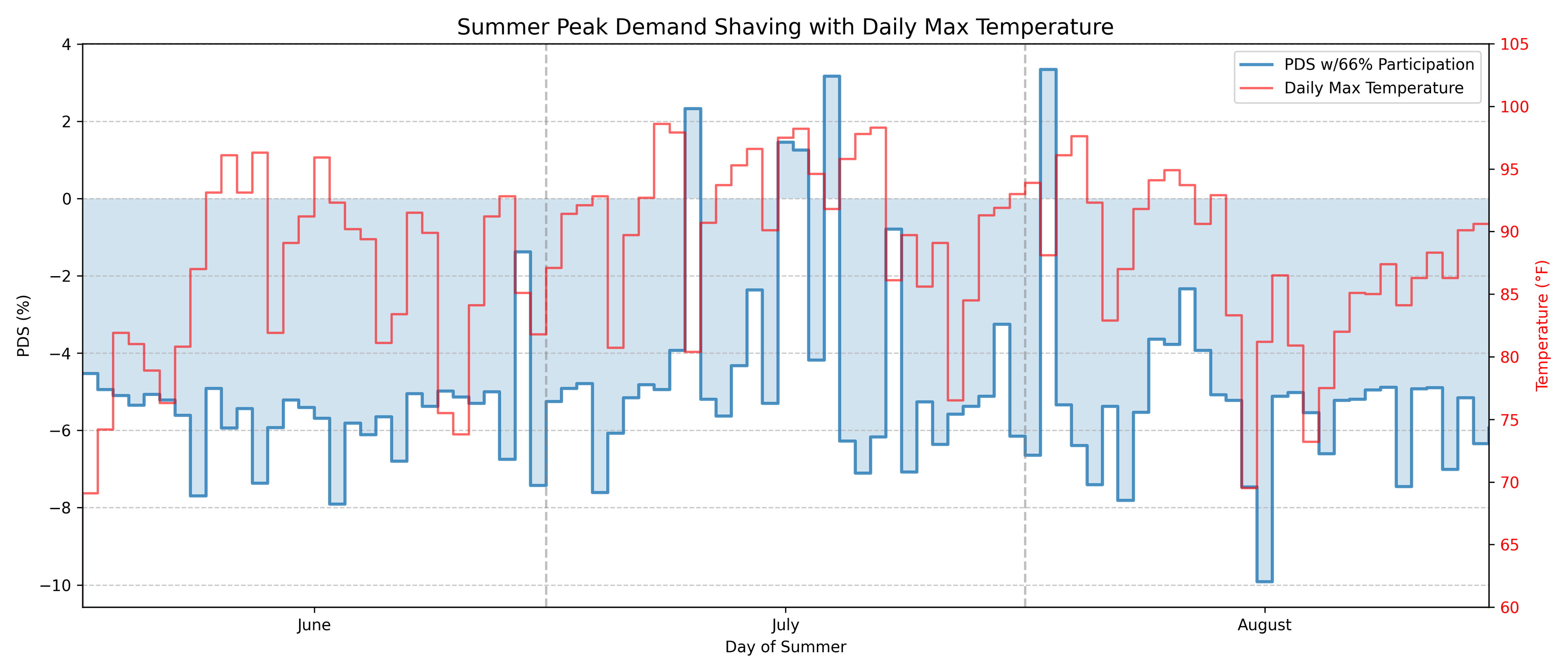}}
\caption{PDS performance loss compare to nominal case over summer months (lower) with HVAC as the only control device, in comparison to the nominal case.}
\label{fig:hvac_PDS}
\end{figure}

The summer results for the HVAC-only scenario are presented in Fig.~\ref{fig:hvac_PDS}. Analysis shows that limiting control to HVAC systems results in a PDS\% loss of 4.75\% compared to the nominal multi-device case. This loss is consistent with the baseline PDS\% observed during non-summer months when HVACs have minimal usage (see Fig.~\ref{fig:nominal_LF}).

Even when participation is restricted solely to HVAC systems—--aligning with some current utility control strategies—--the system still delivers a robust peak demand reduction. Specifically, by subtracting the 4.75\% performance loss from the 16.01\% nominal summer PDS\%, we still achieve a substantial 11--12\% PDS improvement. This indicates that while comprehensive device integration can provide better demand response results, HVAC-focused programs remain a highly effective tool for significant grid-level peak shaving.

\subsubsection{Scalability Analysis}

\begin{figure}[!h]
\centerline{\includegraphics[width=.7\textwidth]{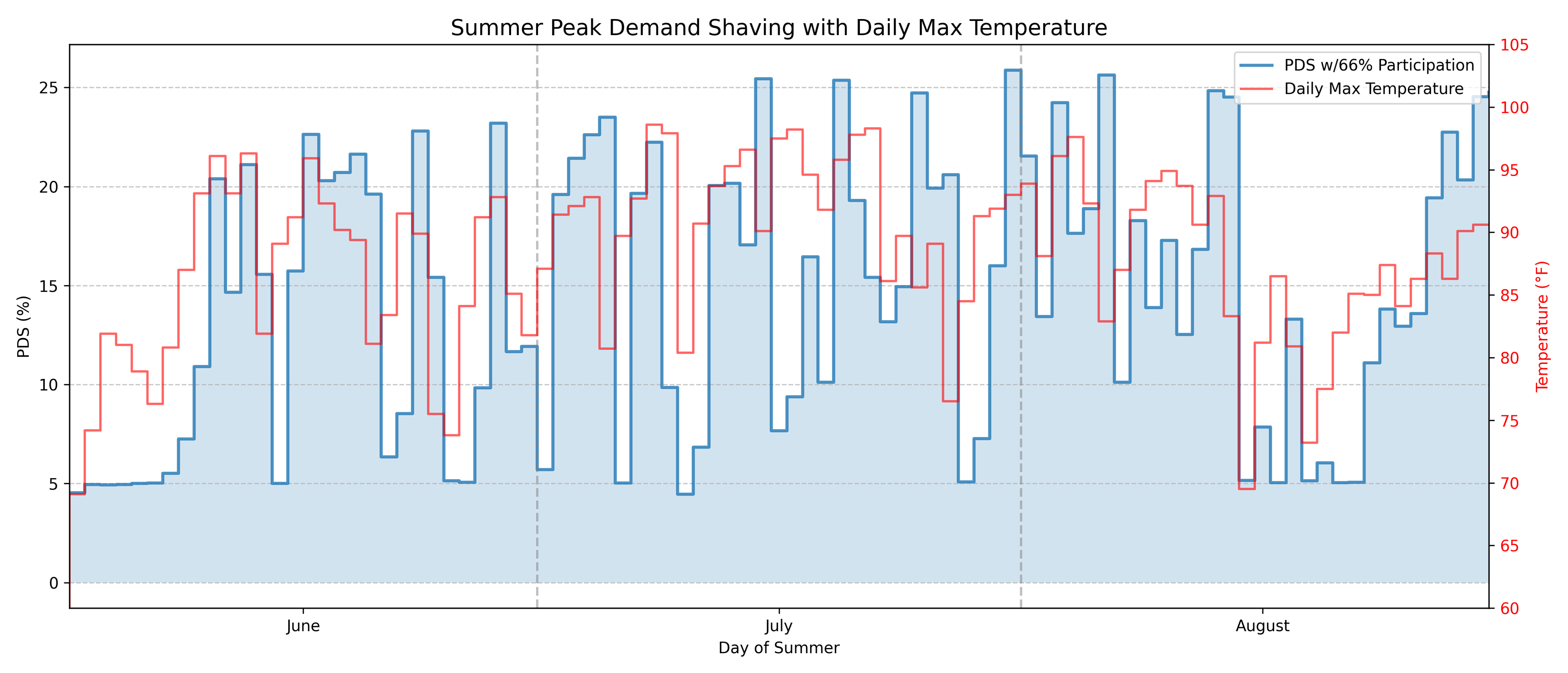}}
\caption{PDS improvement over summer months (lower) as we scale up the size of the system by 30 times.}
\label{fig:scalability}
\end{figure}

To evaluate the computational and performance scalability of the proposed framework, we expanded the network to include 14,070 HEMS units—a 30-fold increase over the nominal case. The resulting summer PDS\% improvements are illustrated in Fig.~\ref{fig:scalability}.
As summarized in Table~\ref{tab:scalability_comparison}, the large-scale simulation yielded a summer PDS\% of 15.96\% and a variation reduction of 19.07\%. These metrics closely align with the nominal case results (16.01\% and 19.33\%, respectively), demonstrating that the algorithm's effectiveness remains consistent even as the system dimensionality increases significantly.

The negligible deviation in performance across different network sizes confirms that the proposed framework is \textbf{highly scalable}. It can effectively manage large-scale residential aggregations without performance degradation, making it suitable for real-world utility applications.

\begin{table}[h]
    \centering
    \caption{Scalability Comparison: Nominal Case vs. 30x Scale-up}\label{tab:scalability_comparison}
    \begin{tabular}{lcc}
        \toprule
        \textbf{Scenario} & \textbf{Summer PDS\%} & \textbf{Variation Reduction \%} \\
        \midrule
        Nominal Case (469 HEMS)  & 16.01 & 19.33 \\
        Large Scale (14,070 HEMS) & 15.96 & 19.07 \\
        \bottomrule
    \end{tabular}
\end{table}

\begin{table}[h!]
    \centering
    \caption{Scenario Analysis Results}
    \label{tab:scenario_analysis}
    \begin{tabular}{lccc}
        \toprule
        \textbf{Scenario} & \textbf{Average PDS} & \textbf{Average variation reduction} & \textbf{AMPS range} \\
        \midrule
        Denver, 66\% participation, 20\% PV/BESS penetration & 16\% & 19\% & 10\%-19\% \\
        Denver, 33\% participation, 20\% PV/BESS penetration & 8\% & 12\% & 5\%-10\% \\
        Denver, 100\% participation, 20\% PV/BESS penetration & 21\% & 20\% & 10\%-19\% \\
        Denver, HVAC only & 12\% & 22\% & 9\%-13\% \\
        Denver, 66\% participation, 60\% PV/BESS penetration & 43\% & 44\% & 27\%-43\% \\
        Denver, 14k customers & 16\% & 19\% & 10\%-16\% \\
        Los Angeles, 66\% participation, 20\% PV/BESS penetration & 17\% & 21\% & 13\%-27\% \\
        Phoenix, 66\% participation, 20\% PV/BESS penetration & 14\% & 40\% & 10\%-20\% \\
        \bottomrule
    \end{tabular}
\end{table}

\section{Conclusion}

The transition toward a highly electrified, DER-rich grid presents both an unprecedented opportunity for flexibility and a critical operational challenge. As demonstrated in this study, harnessing this flexibility using traditional static TOU tariffs inevitably leads to device synchronization, creating secondary ``timer peaks'' that threaten distribution network stability. 

To address this, we proposed and rigorously validated a privacy-preserving, one-way dynamic signaling framework designed to safely unlock deep demand flexibility. By formulating the utility-consumer interaction as a Stackelberg game and employing a feedback-based, context-enriched soft-clustering algorithm, our framework dynamically shapes day-ahead price profiles based solely on aggregate substation feedback and environmental forecasts. This approach effectively circumvents the prohibitive infrastructure costs, computational overhead, and privacy concerns associated with two-way communication and centralized direct load control architectures.

Through high-fidelity simulations utilizing a realistic SMART-DS distribution network populated with diverse behind-the-meter assets, the proposed mechanism demonstrated exceptional performance. In the nominal summer scenario, the dynamic pricing framework successfully mitigated device synchronization, achieving a 16.01\% reduction in peak demand and a 19.33\% reduction in maximum hourly load variation. Notably, an isolated evaluation of an HVAC-only scenario confirmed that intelligent space conditioning serves as the dominant driver of this flexibility, providing the vast majority of the load-shifting capability even before other heterogeneous assets are introduced. On average over a simulated year, the soft-clustering framework improves the daily load factor by 18.2\% and reduces the peak-to-average ratio by 15.4\% compared to the uncoordinated baseline. Crucially, total daily energy consumption dropped by only 4.2\%, confirming that the framework primarily drives healthy load \textit{shifting} rather than disruptive load \textit{shedding}, thereby maintaining end-user utility. 

Furthermore, the framework proved to be both highly scalable and near-optimal. The proposed one-way communication architecture reduces computational overhead significantly compared to iterative two-way negotiation models, while still capturing over 85\% of the theoretical peak-shaving benefits. Importantly, a comprehensive \textbf{participation rate study} revealed non-linear returns on adoption; the framework captures the vast majority of its grid-balancing benefits at a moderate 66\% participation rate. This demonstrates that utilities can achieve high-impact network stabilization without the prohibitive marketing and deployment costs required to drive universal program adoption. 

Our study also confirmed the framework's robustness across massive system scales (over 14,000 HEMS) and its adaptability to extreme climatic zones, such as Phoenix, where it significantly smoothed load volatility even when absolute peak reduction was constrained by saturated cooling demands. Finally, the framework exhibited powerful synergies with renewable integration, nearly doubling its peak-shaving efficacy (to 32.79\%) under higher PV and battery penetration scenarios.

\section*{Acknowledgment}
\bibliographystyle{IEEEtran}
\bibliography{lib}

@IEEEtranBSTCTL{IEEEexample:BSTcontrol,
CTLdash_repeated_names = "no",
}

@techreport{renewhome2025,
  title={Scaling Residential Demand Response by Prioritizing Comfort: Evidence from 25 Million Energy Shifts},
  author={Blasnik, Michael and Reid, Daniel and Lanzisera, Steven},
  institution={Renew Home, LLC},
  type={White Paper},
  year={2025},
  month={Nov}
}

@article{mathieu2025new,
  title={A New Definition and Research Agenda for Demand Response in the Distributed Energy Resource Era},
  author={Mathieu, Johanna L. and Verbi{\v{c}}, Gregor and Morstyn, Thomas and others},
  journal={IEEE Transactions on Energy Markets, Policy, and Regulation},
  volume={3},
  number={3},
  pages={324--339},
  year={2025},
  publisher={IEEE}
}

@article{vardakas2015survey,
  title={A survey on demand response programs in smart grids: Pricing methods and optimization algorithms},
  author={Vardakas, John S. and Zorba, Nizar and Verikoukis, Christos V.},
  journal={IEEE Communications Surveys \& Tutorials},
  volume={17},
  number={1},
  pages={152--178},
  year={2015},
  publisher={IEEE}
}

@article{samad2016automated,
  title={Automated demand response for smart buildings and microgrids: The state of the practice and research challenges},
  author={Samad, Tariq and Koch, Ed and Stluka, Petr},
  journal={Proceedings of the IEEE},
  volume={104},
  number={4},
  pages={726--744},
  year={2016},
  publisher={IEEE}
}

@article{kok2016society,
  title={A society of devices: Integrating intelligent distributed resources with transactive energy},
  author={Kok, Koen and Widergren, Steve},
  journal={IEEE Power and Energy Magazine},
  volume={14},
  number={3},
  pages={34--45},
  year={2016},
  publisher={IEEE}
}

@inproceedings{li2011optimal,
  title={Optimal demand response based on utility maximization in power networks},
  author={Li, Na and Chen, Lijun and Low, Steven H},
  booktitle={2011 IEEE power and energy society general meeting},
  pages={1--8},
  year={2011},
  organization={IEEE}
}

@article{stai2017dispatching,
  title={Dispatching stochastic heterogeneous resources accounting for grid and battery losses},
  author={Stai, Eleni and Reyes-Chamorro, Lorenzo and Sossan, Fabrizio and Le Boudec, Jean-Yves and Paolone, Mario},
  journal={IEEE Transactions on Smart Grid},
  volume={9},
  number={6},
  pages={6522--6539},
  year={2017},
  publisher={IEEE}
}

@article{parvania2013optimal,
  title={Optimal demand response aggregation in wholesale electricity markets},
  author={Parvania, Masood and Fotuhi-Firuzabad, Mahmud and Shahidehpour, Mohammad},
  journal={IEEE Transactions on Smart Grid},
  volume={4},
  number={4},
  pages={1957--1965},
  year={2013},
  publisher={IEEE}
}

@article{nazir2022grid,
  title={Grid-aware aggregation and realtime disaggregation of distributed energy resources in radial networks},
  author={Nazir, Nawaf and Almassalkhi, Mads},
  journal={IEEE Transactions on Power Systems},
  volume={37},
  number={3},
  pages={1706--1717},
  year={2022},
  publisher={IEEE}
}

@article{granitsas2025controlling,
  title={Controlling air conditioners for frequency regulation: A real-world example},
  author={Granitsas, I. M. and others},
  journal={IEEE Transactions on Smart Grid},
  volume={16},
  number={2},
  pages={1221--1232},
  year={2025},
  publisher={IEEE}
}

@article{morstyn2019bilateral,
  title={Bilateral contract networks for peer-to-peer energy trading},
  author={Morstyn, Thomas and Teytelboym, Alexander and McCulloch, Malcolm D.},
  journal={IEEE Transactions on Smart Grid},
  volume={10},
  number={2},
  pages={2026--2035},
  year={2019},
  publisher={IEEE}
}

@article{kim2020p2p,
  title={A P2P-dominant distribution system architecture},
  author={Kim, J. and Dvorkin, Y.},
  journal={IEEE Transactions on Power Systems},
  volume={35},
  number={4},
  pages={2716--2725},
  year={2020},
  publisher={IEEE}
}

@article{molzahn2017survey,
  title={A survey of distributed optimization and control algorithms for electric power systems},
  author={Molzahn, Daniel K. and others},
  journal={IEEE Transactions on Smart Grid},
  volume={8},
  number={6},
  pages={2941--2962},
  year={2017},
  publisher={IEEE}
}

@techreport{zandi2018home,
  title={Home energy management systems: An overview},
  author={Zandi, Helia and Kuruganti, Teja and Vineyard, Edward A. and Fugate, David},
  institution={Oak Ridge National Lab. (ORNL)},
  year={2018}
}

@misc{zippenfenig2023open,
  title={Open-Meteo.com Weather API},
  author={Zippenfenig, Patrick},
  year={2023},
  publisher={Zenodo}
}

@misc{xcel2023residential,
  title={Residential Time of Use Rate Design},
  author={{Xcel Energy}},
  year={2023},
  howpublished={Colorado},
  note={[Online]. Available: https://www.xcelenergy.com}
}

@techreport{nrel2023smartds,
  title={SMART-DS: Synthetic Models for Advanced, Realistic Testing: Distribution Systems and Scenarios},
  author={{National Renewable Energy Laboratory}},
  institution={National Renewable Energy Laboratory},
  year={2023},
  note={[Online]. Available: https://www.nrel.gov/grid/smart-ds.html}
}

@article{parikh2014proximal,
  title={Proximal algorithms},
  author={Parikh, Neal and Boyd, Stephen and others},
  journal={Foundations and Trends{\textregistered} in Optimization},
  volume={1},
  number={3},
  pages={127--239},
  year={2014},
  publisher={Now Publishers, Inc.}
}

@article{zhou2019online,
  title={Online stochastic optimization of networked distributed energy resources},
  author={Zhou, Xinyang and Dall'Anese, Emiliano and Chen, Lijun},
  journal={IEEE Transactions on Automatic Control},
  volume={65},
  number={6},
  pages={2387--2401},
  year={2019},
  publisher={IEEE}
}

@ARTICLE{bernstein2019,
  author={Bernstein, Andrey and Dall’Anese, Emiliano},
  journal={IEEE Transactions on Control of Network Systems}, 
  title={Real-Time Feedback-Based Optimization of Distribution Grids: A Unified Approach}, 
  year={2019},
  volume={6},
  number={3},
  pages={1197-1209},
  doi={10.1109/TCNS.2019.2929648}}

@misc{mehrabi2024optimalmechanismsdemandresponse,
      title={Optimal Mechanisms for Demand Response: An Indifference Set Approach}, 
      author={Mohammad Mehrabi and Omer Karaduman and Stefan Wager},
      year={2024},
      eprint={2409.07655},
      archivePrefix={arXiv},
      primaryClass={math.OC},
      url={https://arxiv.org/abs/2409.07655}, 
}

@inproceedings{EfficientProjections,
author = {Duchi, John and Shalev-Shwartz, Shai and Singer, Yoram and Chandra, Tushar},
title = {Efficient projections onto the l1-ball for learning in high dimensions},
year = {2008},
isbn = {9781605582054},
publisher = {Association for Computing Machinery},
address = {New York, NY, USA},
url = {https://doi.org/10.1145/1390156.1390191},
doi = {10.1145/1390156.1390191},
booktitle = {Proceedings of the 25th International Conference on Machine Learning},
pages = {272–279},
numpages = {8},
location = {Helsinki, Finland},
series = {ICML '08}
}

\appendices   

\section{Algorithm Description}
In this appended section, we will go over details of \cref{alg1} and \cref{alg2}, including objective functions, constraints, training methods and tuning parameters.

\subsection{Objective and Constraints}

Mehrabi et al.\cite{mehrabi2024optimalmechanismsdemandresponse} studied the dual characterization of optimal price signal under the utility-household Stackelberg game mentioned above. Let $\bm{\alpha}$ be the hourly energy price, $\bm{p}_i$ be hourly power consumption of the $i^{th}$ household. The two participants solve the following optimization problems respectively:

\begin{align}
\label{eq:grid}
    &\text{Grid optimization: }\min_{\bm{a}} \rho\left(\sum _{i=1}^n \bm{p}_i\right)\\
\label{eq:HEMS}
    &\text{HEMS optimization: }\min_{\bm{p}_i} \bm{a}^T \bm{p}_i\quad s.t.\,\,\bm{p}_i\in \Omega_i
\end{align} 
where the grid cost $\rho(\cdot)$ is an order-one homogeneous function of the total household consumption, and $\Omega_i$ is a pre-specified feasible set by the households.

Let $h(\bm{\alpha}; \Omega_i) = \min\{\bm{a}^T \bm{p}_i: \bm{p}_i\in \Omega_i\}$ be the optimum of \eqref{eq:grid}, then Mehrabi et al.\cite{mehrabi2024optimalmechanismsdemandresponse} showed that the mean-field solution to \eqref{eq:HEMS} can be characterized as follows:
\begin{equation} \label{eq:basic}
    \bm{\alpha}^* = \arg\max_{\bm{\alpha} \in \mathcal{A}} \left\{H(\bm{\alpha})\right\}
\end{equation}
where $\mathcal{A}=\left\{\bm{z}\in \mathbb{R}^d :\bm{z}^T \bm{p} \leq \rho(\bm{p}), \, \forall \bm{p} \in  \mathbb{R}^d \right\}$, and $H(\bm{\alpha})=\mathbb{E}[h(\bm{\alpha}; \Omega_i)]$. In our case, $d=48$ since the system adopts two-day-ahead planing.

Since it is not stated explicitly in the proceeding work, we give the following lemma showing that the constraint $\bm{\alpha}\in\mathcal{A}$ is equivalent to projection onto its dual-norm ball:

\begin{lemma} \label{lemma:dualnorm}
With the grid cost function $\rho$ characterized by \eqref{eq:grid}, the constraint set is equivalent to:

\begin{equation}
   \mathcal{A}=\{\bm{z}\in \mathbb{R}^d :\rho^*(\bm{z}) \leq 1\}
\end{equation}
where $\rho^*$ is the dual norm of $\rho$.
\end{lemma}

We give a few examples to demonstrate how to simplify the constraints using dual norm:
\begin{example}[$l_\infty$ Norm]\label{ex:infnorm}
Suppose $\rho(\bm{p}) = ||\bm{p}||_\infty$, then its dual norm $\rho^*(\bm{z}) = ||\bm{z}||_1$. We can then utilize any known efficient algorithms to project output onto the $l_1$-ball, such as in Duchi et al.\cite{EfficientProjections}. 
\end{example}

\begin{example}[$l_2$ Norm with Smoothness Regularization]
As discussed in \cref{utility_objective}, a utility may have multiple objectives, including peak hour demand and smoothness of the entire
demand profile. Suppose $\rho(\bm{p}) = (||\bm{p}||_2^2 + \lambda \bm{p}^TD\bm{p})^{1/2} :=||\bm{p}||_K$, where $D \in \mathbb{R}^{d \times d}$ is a Toeplitz matrix taking the difference of two neighboring entries, i.e. $\bm{p}^TD\bm{p} = (\bm{p}_{d}-\bm{p}_1)^2 + \sum_{j=1}^{d-1} (\bm{p}_{j+1}-\bm{p}_j)^2$, $K := I + \lambda D^TD$. Such cost function penalizes both hourly changes in $\bm{p}$ and the magnitude of $\bm{p}$. By \cref{lemma:dualnorm}, $\mathcal{A} =\{\bm{z} \in \mathbb{R}^{d}:\bm{z}^T K^{-1}\bm{z}\leq1\}$. Projection onto such set can be formally written as an optimization problem:
\begin{equation*}
    \arg\min_{\bm{\alpha}} ||\bm{z}-\bm{\alpha}||_2\qquad s.t.\ \bm{\alpha}\in\mathcal{A}
\end{equation*}
which can be solved by bisection method:
\begin{equation*}
    \bm{\alpha}^* = 
    \begin{cases}
    \bm{z},\quad \textnormal{ if } \bm{z}^TK^{-1}\bm{z} \leq 1\\ (I+\mu K^{-1})^{-1}\bm{z},
    \quad \textnormal{ if } \bm{z}^TK^{-1}\bm{z} > 1
    \end{cases}
\end{equation*}
where $(I+\mu K^{-1})^{-1}$ is the direct application of first order condition on the Lagrange function, and we need to solve for $\mu>0,\; s.t.\; \bm{\alpha}^{*T}K^{-1}\bm{\alpha}^{*} =1$ in the second case. This particular choice of grid cost function is used in our simulation.
\end{example}

Now suppose weather forecasts $X_t\in \mathcal{X}$ are available, we define $\bm{\alpha}_\theta: \mathcal{X} \mapsto \mathbb{R}^d$ to be the pricing function parametrized by $\theta$, i.e.
\begin{equation}
\label{eq:pricing_func}
    \bm{\alpha}_t = \bm{\alpha}_\theta(X_t).
\end{equation}

The following proposition leverages the mean-field solution \eqref{eq:basic} to solve for the optimal price signal of contextual learning:
\begin{proposition}
\label{prop:sol}
    Under assumption \eqref{eq:pricing_func} on parametrized pricing function, the optimal price of contextual learning can be characterized as:
\begin{equation}
\label{eq:optimalsol}
    \theta^* = \arg\max_{\theta \in \Theta} \{H(\bm{\alpha}_\theta(X_t))\}
\end{equation}
where $\Theta$ is the feasible set of parameters such that the pricing function $\bm{\alpha}_\theta$ maps $\mathcal{X}$ to $\mathcal{A}$.
\end{proposition}

In practice, it is often challenging to convert constraints from $\bm{\alpha} \in \mathcal{A}$ to $\theta \in \Theta$, hence an explicit characterization of the feasible set $\Theta$ can be intractable. To address this, we propose an algorithm for learning tasks with constrained predictions in the next section.

\subsection{Methodology and Algorithms}

Existing literature has been focusing on scenarios where constraints on the predictions can be converted to constraints on the training variables, then utilizing tools like projected gradient descent to get convergence guarantee \cite{parikh2014proximal}. However, in fact it is hard to obtain a simple closed-form constraint set of the training variables in many cases, such as the weight matrix in a neural network where prediction mapping is not bijective. In those cases, simply projecting outputs onto the constraint set is not guaranteed to reach convergence. The issue is prevalent even when the loss is convex, as the composition of loss and projection usually does not preserve convexity.

In this section, we propose a soft-clustering neural network architecture as a general solution to contextual learning with prediction constraints. It has two main components: first is to train a neural network that takes in forecast data to predict the distribution of a future day type; second is to use realized hourly demand from power system outputs to optimize the price profile corresponding to each cluster. The final algorithm output is a weighted average of cluster price profiles. 

Such architecture avoids the conversion of constraint set from $\mathcal{A}$ to $\Theta$, and only calls for projected gradient descent to impose the constraints on each cluster. Formally, we define the soft-clustering neural network architecture:

\begin{definition} \label{nnstruc}
Assume $\psi(X_t)$ is a fully connected neural network parametrized by $\theta = (w_1,\, b_1,\, w_2,\, b_2)$, then the soft-clustering neural network is defined as follows:

\begin{equation}
\begin{aligned}
    \label{eq:softclass}
     \bm{\alpha}_t &= \bm{\alpha}_\theta(X_t) = \kappa \cdot \psi(X_t)\\ \psi(X_t) &= \text{softmax}\left\{\left(w_2h\left(w_1X_t+b_1\right)\right)+b_2\right\}
\end{aligned}
\end{equation}
where the clusters $\kappa \in \mathbb{R}^{d\times k}$, the neural network output $\psi(X)\in \mathbb{R}^k$, and the final price prediction $\bm{\alpha}_t \in \mathcal{A}\subset\mathbb{R}^d$. 
\end{definition}

Candidates for the activation function $h$ include sigmoid, ReLU, hyperbolic tangent, etc. In our simulation, we choose to use hyperbolic tangent $\tanh$ as the activation function.

The following lemma establishes the equivalency of constraints on final predictions and on cluster predictions, i.e. every column in $\kappa$ needs to satisfy the constraints. It plays a key role in addressing the challenge of constrained prediction problem:

\begin{lemma}
Under the neural network structure in \eqref{eq:softclass}, the following statements are equivalent:
    \begin{equation*}
    \begin{aligned}
    \bm{\alpha}_t &= \kappa \cdot \psi(X_t) \in \mathcal{A}\\ &\Leftrightarrow \kappa \in \{\kappa| \kappa \cdot v \in \mathcal{A},\, \forall v \in \mathbb{S}^k\}\\ &\Leftrightarrow \kappa \in \{\kappa| \kappa_{\cdot j} \in \mathcal{A}, \forall 1\leq j \leq k\}
    \end{aligned}
    \end{equation*}
where $\mathbb{S}^k$ is a $k$-dimensional simplex.
\end{lemma}

An advantage of the soft-clustering architecture is training separation. We update the model $\psi(x)$ and the clusters $\kappa$ via different methods: parameters $(w_1, w_2, b_1, b_2)$ from $\psi(x)$ are trained by unconstrained gradient descent offline, while constrained clusters $\kappa$ are updated by projected gradient descent (PGD) online. The loss function for each component is defined as follows:

For offline clustering: 
\begin{equation*}
\begin{aligned}
&\mathcal{L}^{off} = \frac{1}{n}\sum_{i=1}^n\left[||\hat{y}_i^{(temp)}-y_i^{(temp)}||^2+||\hat{y}_i^{(solar)}-y_i^{(solar)}||^2\right]\\
    &\hat{y}_i^{(temp)} = \mu^{(temp)}\psi(X_i) = \sum_{k=1}^{K}\psi(X_i)_k\mu_{k}^{(temp)}\\
    &\hat{y}_i^{(solar)} = \mu^{(solar)}\psi(X_i)=\sum_{k=1}^{K}\psi(X_i)_k\mu_{k}^{(solar)}\\
\end{aligned}
\end{equation*}
where the features $y^{(temp)}$ and $y^{(solar)}$ are normalized; $\mu^{(temp)}$ and $\mu^{(solar)}$ are matrices of cluster centroids, for temperature and solar irradiance respectively, and are learned simultaneously during the neural network training. To further encourage divergence and sensitivity of clustering weights, we also add optional entropy loss and contrast loss to $\mathcal{L}^{off}$:

\begin{equation*}
\begin{aligned}
\mathcal{L}_{\text{entropy}}
&= - \frac{1}{n} \sum_{i=1}^{n} \sum_{j=1}^k 
w_{ij}\,\log w_{ij} \\
\mathcal{L}_{\text{contrast}}
&= - \frac{\sum_{i=1}^{n}\sum_{j=1}^{n}
\left\| w_i - w_j \right\|_2^2}
{\sum_{i=1}^{n}\sum_{j=1}^{n}
\left\| x_i - x_j \right\|_2^2} .
\end{aligned}
\end{equation*}

For online update: 
\begin{equation*}
    \begin{aligned}
    \max_{\kappa_{\cdot j} \in \mathcal{A},\; \forall j} \left\{H(\kappa \cdot \psi(X_t))\right\}
    \end{aligned}
\end{equation*}
where $\mathcal{A} =\{z \in \mathbb{R}^{d}:z^T K^{-1}z\leq1\}$, $K = I + \lambda D^TD$. We leverage the first-order algorithm proposed by Mehrabi et al.\cite{mehrabi2024optimalmechanismsdemandresponse} assuming homogeneous day types to update each cluster, where the gradient of cluster $j$ at each step $t$ is:

\begin{equation*}
    g_t^j = \psi_j(X_t)\left(\frac{1}{n}\sum_{i=1}^n \bm{p}_i(\bm{\alpha}_t, X_t)\right)
\end{equation*}

\cref{alg1} gives pseudo-code of the soft-clustering architecture. Note that we utilize projected gradient descent (PSD) in our demonstration, future improvements include adding (adaptive) momentum term to reduce sensitivity to tuning and stabilize the training.

\begin{center}
\begin{minipage}{0.75\textwidth}
\begin{algorithm}[H]
\DontPrintSemicolon
\setlength{\algomargin}{0pt} 
\setlength{\algowidth}{0.75\textwidth} 
\newcommand{\rcomment}[1]{\hfill\parbox[t]{3cm}{\raggedleft\footnotesize\textit{#1}}}

\caption{Soft-Clustering Neural Network Algorithm} \label{alg3}

$w, b, \kappa \gets w_0, b_0, \kappa_0$\tcp*{Initialize training variables}

\For{$s \gets 1$ \KwTo $S$}{
    $w_s \gets w_{s-1} - \gamma_s \cdot \nabla_w \mathcal{L}^{off}(w, b_{s-1}, \mu_{s-1})$
    
    $b_s \gets b_{s-1} - \gamma_s \cdot \nabla_b \mathcal{L}^{off}(w_{s-1}, b, \mu_{s-1})$
    
    $\mu_s \gets \mu_{s-1} - \gamma_s \cdot \nabla_\mu \mathcal{L}^{off}(w_{s-1}, b_{s-1}, \mu)$
}
\tcp*{Train model $\psi(x)$ without constraints}

\For{$t \gets 1$ \KwTo $T$}{
    $\kappa_t' \gets \kappa_{t-1} + \eta_t \cdot g_t$
    
    $\kappa_t \gets \mathrm{Proj}_\mathcal{A}(\kappa_t')$
}
\tcp*{Train clusters $\kappa$ using PGD}
\end{algorithm}
\end{minipage}
\end{center}

\cref{tab:summary} summarizes all tuning parameters in the soft-clustering algorithm and the settings in our simulation. We pick $\eta_t = \eta/||g_t||_2$, i.e. normalize the step size by $l_2$-norm of the gradient. Such normalization preserves the relative magnitude of within-day demand, and significantly reduces the sensitivity of convergence on step size. The remaining hyperparameters are selected based on empirical calibration. Across a broad set of test cases, we find that the chosen configuration delivers stable and consistently strong performance, suggesting that the method is not unduly sensitive to moderate variations.

\vspace{3mm}
\begin{table}[H]
\begin{center}
\caption{Tuning Parameters of Soft-Clustering Algorithm}
\label{tab:summary}
\renewcommand{\arraystretch}{1}
\begin{tabular}{|L{2.7cm}|Y{4.3cm}|R{2.3cm}|}
 \hline
 \textbf{Parameter} & \textbf{Description}& \textbf{Simulation}\\
 \hline
 Number of clusters $k$ &  Reflects expectation on \# day types& $k=6$\\ 
 \hline
 Hidden layer nodes $w_2$& Dimension of hidden layer in neural network $\psi$ & $d_{in} = d_{out} = 48$, $d_{hidden} = 60$\\ 
 \hline
 Online penalty $\lambda$&  Controls penalty on variation in grid cost& 
 $\lambda_{l_2\text{-}norm} = 0.1$, $\lambda_{variation} = 0.9$\\ 
 \hline
 Offline penalty  $\lambda$&  Controls divergence and sensitivity of clustering weights& $\lambda_{entropy} = -0.4$, $\lambda_{contrast} = 0.5$\\ 
 \hline
 Learning rate $\gamma_s$, $\eta_t$ & Step size of offline and online training & $\gamma_s=0.001$, $\eta_t = 0.1/||g_t||_2$\\
 \hline
\end{tabular}
\end{center}
\end{table}

\end{document}